\font \ninerm=cmr9

\font \eightsl=cmsl8

\font \bigbf=cmbx10 scaled \magstep1
\font \Bigbf=cmbx10 scaled \magstep2

\font \tengoth=eufm10
\font \sevengoth=eufm7
\font \fivegoth=eufm5

\newfam\gothfam
\textfont \gothfam=\tengoth
\scriptfont \gothfam=\sevengoth
\scriptscriptfont \gothfam=\fivegoth

%
\font \tenmath=msbm10
\font \sevenmath=msbm7
\font \fivemath=msbm5

\newfam\mathfam
\textfont \mathfam=\tenmath
\scriptfont \mathfam=\sevenmath
\scriptscriptfont \mathfam=\fivemath
\def\math{\fam\mathfam\tenmath}

%
%
\def\titre#1{\centerline{\Bigbf #1}\nobreak\nobreak\vglue 10mm\nobreak}

\def\paragraphe#1{\bigskip\goodbreak {\bigbf #1}\nobreak\vglue 12pt\nobreak}
\def\alinea#1{\medskip\allowbreak{\bf#1}\nobreak\vglue 9pt\nobreak}
\def\ssq{\smallskip\qquad}
\def\msq{\medskip\qquad}

%
%
\def\th#1{\bigskip\goodbreak {\bf Th\'eor\`eme #1.} \par\nobreak \sl }
\def\prop#1{\bigskip\goodbreak {\bf Proposition #1.} \par\nobreak \sl }
\def\lemme#1{\bigskip\goodbreak {\bf Lemme #1.} \par\nobreak \sl }
\def\cor#1{\bigskip\goodbreak {\bf Corollaire #1.} \par\nobreak \sl }
\def\definition#1{\bigskip\goodbreak {\bf D\'efinition #1.} \par\nobreak \sl }
\def\dem{\bigskip\goodbreak \it D\'emonstration. \rm}
\def\ndem{\bigskip\goodbreak \rm}
\def\qed{\par\nobreak\hfill $\bullet$ \par\goodbreak}
%
%
\def\uple#1#2{#1_1,\ldots ,{#1}_{#2}}
\def \restr#1{\mathstrut_{\textstyle |}\raise-6pt\hbox{$\scriptstyle #1$}}
\def \srestr#1{\mathstrut_{\scriptstyle |}\hbox to -1.5pt{}\raise-4pt\hbox{$\scriptscriptstyle #1$}}
\def \inver{^{-1}}
\def\dbar{d\!\!\hbox to 4.5pt{\hfill\vrule height 5.5pt depth -5.3pt
        width 3.5pt}}

\def\frac#1#2{{\textstyle {#1\over #2}}}
\def\R{{\math R}}
\def\C{{\math C}}

\def\N{{\math N}}
\def\Z{{\math Z}}

\let\wt=\widetilde
\def\permuc#1#2#3#4{#1#2#3#4+#1#3#4#2+#1#4#2#3}

\def\fleche#1{\mathop{\hbox to #1 mm{\rightarrowfill}}\limits}
\def\gfleche#1{\mathop{\hbox to #1 mm{\leftarrowfill}}\limits}
\def\inj#1{\mathop{\hbox to #1 mm{$\lhook\joinrel$\rightarrowfill}}\limits}
\def\ginj#1{\mathop{\hbox to #1 mm{\leftarrowfill$\joinrel\rhook$}}\limits}
\def\surj#1{\mathop{\hbox to #1 mm{\rightarrowfill\hskip 2pt\llap{$\rightarrow$}}}\limits}
\def\gsurj#1{\mathop{\hbox to #1 mm{\rlap{$\leftarrow$}\hskip 2pt \leftarrowfill}}\limits}
%
%
\def \g#1{\hbox{\tengoth #1}}
\def \sg#1{\hbox{\sevengoth #1}}

\def\Cal #1{{\cal #1}}
%
%

\def \mop#1{\mathop{\hbox{\rm #1}}\nolimits}

\def \mopl#1{\mathop{\hbox{\rm #1}}\limits}

%
%
\def \bib #1{\null\medskip \strut\llap{[#1]\quad}}
\def\cite#1{[#1]}
\input epsf
\long\def\dessin#1#2{\null
               \bigskip
                \begingroup \epsfysize = #1 $$\epsfbox {#2}$$ \endgroup
                \bigskip
                \goodbreak}                
\magnification1000 
\tolerance 1000

\font\hb=cmbx12

\font\bbfnt=msbm10
\def\BbbR{\hbox{\bbfnt\char'122}}

\def\BbbC{\hbox{\bbfnt\char'103}}

\def\BbbP{\hbox{\bbfnt\char'120}}
\def\BbbT{\hbox{\bbfnt\char'124}}

\def\pt{\hbox{\hb .}}

\let\frak=\goth

\parindent=0cm

\null
\titre{Choix des signes pour la formalit\'e}
\vskip -6mm
\titre{de M. Kontsevich}
\centerline{D. Arnal
        \footnote{*}{\ninerm Universit\'e de Metz, D\'epartement de Math\'ematiques, \^\i le du Saulcy, 57045 Metz CEDEX 01. arnal@pon\-celet.univ-metz.fr}, 
        D. Manchon 
        \footnote{**}{\ninerm Institut Elie Cartan, CNRS, BP 239, 54506 Vand\oe uvre CEDEX. manchon@iecn.u-nancy.fr}et M. Masmoudi
        \footnote{***}{\ninerm Universit\'e de Metz, D\'epartement de Math\'ematiques, \^\i le du Saulcy, 57045 Metz CEDEX 01. masmoudi@pon\-celet.univ-metz.fr}}

\vskip 4cm
{\baselineskip=11pt
\bf R\'esum\'e~: \rm L'expression explicite de la formalit\'e de M. Kontsevich sur $\R^d$ est la base de la preuve du th\'eor\`eme de formalit\'e pour une vari\'et\'e quelconque \cite {K1 \S~7}, qui implique \`a son tour l'existence d'\'etoile-produits sur une vari\'et\'e de Poisson quelconque. Nous proposons ici un choix coh\'erent d'orientations et de signes qui permet de reprendre la d\'emonstration du th\'eor\`eme pour $\R^d$ en tenant compte des signes qui apparaissent devant les diff\'erents termes de l'\'equation de formalit\'e.
\par}
\paragraphe{Introduction}
\qquad La conjecture de formalit\'e a \'et\'e introduite par M. Kontsevich \cite {K2}~: elle affirme l'existence d'un $L_\infty$-quasi-isomorphisme de $\g g_1$ vers $\g g_2$, o\`u $\g g_1$ et $\g g_2$ sont les deux alg\`ebres de Lie diff\'erentielles gradu\'ees naturellement associ\'ees \`a une vari\'et\'e $M$~: pr\'ecis\'ement $\g g_1$ est l'alg\`ebre de Lie diff\'erentielle gradu\'ee des multi-champs de vecteurs munie de la diff\'erentielle nulle et du crochet de Schouten, et $\g g_2$ est l'alg\`ebre de Lie diff\'erentielle gradu\'ee des op\'erateurs polydiff\'erentiels munie de la diff\'erentielle de Hochschild et du crochet de Gerstenhaber. 
\ssq
Les \'el\'ements de degr\'e $n$ dans $\g g_1$ sont les $(n+1)$-champs de vecteurs, et les \'el\'ements de degr\'e $n$ dans $\g g_2$ sont les op\'erateurs $(n+1)$-diff\'erentiels. Dans les espaces gradu\'es d\'ecal\'es $\g g_1[1]$ et $\g g_2[1]$ ce sont les $(n+2)$-champs de vecteurs (resp. les op\'erateurs $(n+2)$-diff\'erentiels) qui sont de degr\'e $n$.
\ssq
Toute alg\`ebre de Lie diff\'erentielle gradu\'ee est une $L_\infty$-alg\`ebre. Cela signifie en particulier que les structures d'alg\`ebres de Lie diff\'erentielles gradu\'ees sur $\g g_1$ et $\g g_2$ induisent des cod\'erivations $Q$ et $Q'$ de degr\'e $1$ sur des cog\`ebres $\Cal C(\g g_1)=S^+(\g g_1[1])$ et $\Cal C(\g g_2)=S^+(\g g_2[1])$ respectivement (cf. \S~II.3 et II.4), v\'erifiant toutes deux l'\'equation ma\^\i tresse~:
$$[Q,Q]=0, \hbox to 12 mm{}[Q',Q']=0.$$
Un $L_\infty$-quasi-isomorphisme de $\g g_1$ vers $\g g_2$ est par d\'efinition un morphisme de cog\`ebres~:
$$\Cal U:\Cal C(\g g_1)\longrightarrow \Cal C(\g g_2)$$
de degr\'e z\'ero et commutant aux cod\'erivations, c'est-\`a-dire v\'erifiant l'\'equation~:
$$\Cal U\circ Q=Q'\circ \Cal U,$$
et dont la restriction \`a $\g g_1$ est un quasi-isomorphisme de complexes de $\g g_1$ dans $\g g_2$. M. Kontsevich d\'emontre dans \cite {K1} la conjecture de formalit\'e, c'est-\`a-dire l'existence d'un $L_\infty$-quasi-isomorphisme de $\g g_1$ vers $\g g_2$, pour toute vari\'et\'e $M$ de classe $C^\infty$. La premi\`ere \'etape de la preuve (et m\^eme l'essentiel du travail) consiste en la construction explicite du $L_\infty$-quasi-isomorphisme $\Cal U$ pour $M=\R^d$. Le $L_\infty$-quasi-isomorphisme $\Cal U$ est uniquement d\'etermin\'e par ses coefficients de Taylor~:
$$\Cal U_n:S^n(\g g_1[1])\longrightarrow \g g_2[1].$$
(cf. \S~ III.2). Si les $\alpha_k$ sont des $s_k$-champs de vecteurs, ils sont de degr\'e $s_k-2$ dans l'espace d\'ecal\'e $\g g_1[1]$, et donc $\Cal U_n(\alpha_1\cdots \alpha_n)$ est d'ordre $s_1+\cdots +s_n-2n$ dans $\g g_2[1]$. C'est donc un op\'erateur $m$-diff\'erentiel, avec~:
$$\sum_{k=1}^n s_k=2n+m-2.\eqno{(*)}$$
Les coefficients de Taylor sont construits \`a l'aide de poids et de graphes~: on d\'esigne par $G_{n,m}$ l'ensemble des graphes \'etiquet\'es et orient\'es ayant $n$ sommets du premier type (sommets a\'eriens) et $m$ sommets du deuxi\`eme type (sommets terrestres) tels que~:
\smallskip
1). Les ar\^etes partent toutes des sommets a\'eriens.
\smallskip
2). Le but d'une ar\^ete est diff\'erent de sa source (il n'y a pas de boucles).
\smallskip
3). Il n'y a pas d'ar\^etes multiples.
\ssq
A tout graphe $\Gamma\in G_{n,m}$ muni d'un ordre sur l'ensemble de ses ar\^etes, et \`a tout $n$-uple de multi-champs de vecteurs $\uple \alpha n$ on peut associer de mani\`ere naturelle un op\'erateur $m$-diff\'erentiel $B_\Gamma (\alpha_1\otimes\cdots\otimes \alpha_n)$ lorsque pour tout $j\in\{1,\ldots ,n\}$, $\alpha_j$ est un $s_j$-champ de vecteurs, o\`u $s_j$ d\'esigne le nombre d'ar\^etes qui partent du sommet a\'erien num\'ero $j$ \cite {K1 \S~ 6.3}.
\ssq
Le coefficient de Taylor $\Cal U_n$ est alors donn\'e par la formule~:
$$\Cal U_n(\alpha_1\cdots \alpha_n)=\sum_{\Gamma\in G_{n,m}}
        W_\Gamma B_\Gamma(\alpha_1\otimes\cdots\otimes\alpha_n),$$
o\`u l'entier $m$ est reli\'e \`a $n$ et aux $\alpha_j$ par la formule (*) ci-dessus.
\ssq
Le poids $W_\Gamma$ est nul sauf si le nombre d'ar\^etes $|E_\Gamma|$ du graphe $\Gamma$ est pr\'ecis\'ement \'egal \`a $2n+m-2$. Il s'obtient en int\'egrant une forme ferm\'ee $\omega_\Gamma$ de degr\'e $|E_\Gamma|$ sur une composante connexe de la compactification de Fulton-McPherson d'un espace de configuration $C^+_{A,B}$, qui est pr\'ecis\'ement de dimension $2n+m-2$ \cite {FM}, \cite {K1 \S~5}. Il d\'epend lui aussi d'un ordre sur l'ensemble des ar\^etes, mais le produit $W_\Gamma.B_\Gamma$ n'en d\'epend plus.
\ssq
Pour prouver le th\'eor\`eme de formalit\'e M. Kontsevich montre que le morphisme de cog\`ebres $\Cal U$ dont les coefficients de Taylor sont les $\Cal U_n$ d\'efinis ci-dessus est un $L_\infty$-quasi-isomorphisme. La m\'ethode consiste \`a ramener l'\'equation de formalit\'e $\Cal U\circ Q=Q'\circ \Cal U$, qui se d\'eveloppe \`a l'aide des coefficients de Taylor de $\Cal U$, $Q$ et $Q'$~:
$$\eqalign{Q'_1{\cal U}_n\left(\alpha_1\pt...\pt \alpha_n\right)+&{1\over 2}\sum_{I\sqcup
J=\{1,...,n\}\atop I,J\neq\emptyset}\pm Q'_2\left({\cal U}_{\vert
I\vert}\left(\alpha_I\right)\pt{\cal U}_{\vert J\vert}\left(\alpha_J\right)\right)=\cr
&=\sum_{k=1}^n \pm {\cal U}_n\left(Q_1(\alpha_k)\pt
\alpha_1\pt...\pt\widehat{\alpha_k}\pt...\pt \alpha_n\right)+\cr
&\hskip 0.5cm{1\over 2}\sum_{k\neq
l}\pm {\cal U}_{n-1}\left(Q_2(\alpha_k\pt
\alpha_l)\pt\alpha_1\pt...\pt
\widehat{\alpha_k}\pt...\pt\widehat{\alpha_l}\pt...\pt\alpha_n\right)\cr}$$
\`a l'application de la formule de Stokes pour les formes $\omega_\Gamma$ sur l'ensemble des faces de codimension $1$ du bord des espaces de configuration.
\msq
Nous proposons dans la premi\`ere partie de ce travail un choix d'orientation des espaces de configuration (ou plus exactement d'une composante connexe de ceux-ci) $C^+_{A,B}$, et un choix coh\'erent d'orientation pour chacune des faces de codimension $1$ du bord de la compactification.
\ssq
Au chapitre II nous explicitons l'isomorphisme $\Phi:S^n(\g g[1])\fleche 8^\sim \Lambda^n(\g g)[n]$ mentionn\'e dans \cite {K1 \S~4.2} pour tout espace vectoriel gradu\'e $\g g$, afin de pr\'eciser le passage du langage des alg\`ebres de Lie diff\'erentielles gradu\'ees et des $L_\infty$-alg\`ebres au langage des $Q$-vari\'et\'es formelles gradu\'ees point\'ees (\cite {AKSZ}, \cite {K1 \S~4.1}).
\ssq
Nous donnons au chapitre III une formule explicite pour un champ de vecteurs sur une vari\'et\'e formelle gradu\'ee point\'ee ou un morphisme de vari\'et\'es formelles gradu\'ees point\'ees en fonction de leurs coefficients de Taylor respectifs. La d\'emonstration est de nature combinatoire et se fait en explicitant la restriction \`a la puissance sym\'etrique $n$-i\`eme par r\'ecurrence sur $n$.
\ssq
Dans le chapitre IV nous exprimons les deux alg\`ebres de Lie diff\'erentielles gradu\'ees qui nous int\'eressent comme $Q$-vari\'et\'es formelles gradu\'ees point\'ees. L'isomorphisme d'espaces gradu\'es $\Phi$ explicit\'e au chapitre II est ici essentiel. Pour la suite nous sommes amen\'es \`a modifier l'alg\`ebre de Lie diff\'erentielle gradu\'ee des multi-champs de vecteurs~: nous utilisons un crochet de Lie gradu\'e $[\ ,\ ]'$ li\'e au crochet de Schouten par la formule~:
$$[x,y]'=-[y,x]_{\hbox{\eightsl Schouten}}.$$
Les deux crochets co\"\i ncident modulo un changement de signe en pr\'esence de deux \'el\'ements impairs. Nous pr\'ecisons au paragraphe IV.4 les signes (du type Quillen) qui apparaissent dans l'\'equation de formalit\'e.
\ssq
Enfin nous montrons au chapitre VI que modulo tous les choix effectu\'es pr\'ec\'edemment $\Cal U$ est bien un $L_\infty$-morphisme.
\ssq
Le chapitre V est assez largement ind\'ependant du reste de l'article bien que directement reli\'e \`a \cite {K1}~: nous y donnons une d\'emonstration d\'etaill\'ee du th\'eor\`eme de quasi-inversion des quasi-isomorphismes donn\'e dans \cite {K1 \S~4.4-4.5}. Enfin nous rappelons en appendice le lien entre formalit\'e et quantification par d\'eformation.
\vskip 8mm plus 2mm
\paragraphe{I. Orientation des espaces de configuration}
\alinea{I.1. Trois choix de param\'etrage}
\definition{\rm (Espaces de configuration)}
Soit ${\cal H}$ le demi-plan de Poincar\'e (${\cal H}=\{z\in\BbbC,~{\frak{Im}}z>0\}$).
Appelons $Conf^+\left(\{z_1,...,z_n\};\{t_1,...,t_m\}\right)$ l'ensemble des nuages de points~:

$$\left\{(z_1,...,z_n;t_1,..,t_m),~\hbox{t. q.}~ z_i\in{\cal H},~t_j\in\BbbR,~z_i\neq
z_{i'}~\hbox{si}~i\neq i',~t_1<...<t_m\right\}$$
et $C^+_{\{p_1,...,p_n\};\{q_1,...,q_m\}}$ le quotient de cette vari\'et\'e sous l'action du groupe 
$G$ de toutes les transformations de la forme~:

$$z_i\mapsto az_i+b,~t_j\mapsto at_j+b\qquad (a>0,~b\in\BbbR).$$
\ndem
La vari\'et\'e $C^+_{\{p_1,...,p_n\};\{q_1,...,q_m\}}$ est donc de dimension $2n+m-2$. Cette
vari\'et\'e est par convention orient\'ee par le passage au quotient de la forme~:

$$\Omega_{\{z_1,...,z_n\};\{t_1,...,t_m\}}=dx_1\wedge dy_1\wedge...\wedge dx_n\wedge
dy_n\wedge dt_1\wedge ... \wedge dt_m$$
o\`u $z_j=x_j+iy_j$. Le groupe des transformations consid\'er\'ees pr\'eserve l'orientation. On en
d\'eduit une orientation des espaces $C^+_{\{p_1,...,p_n\};\{q_1,...,q_m\}}$. Plus
pr\'ecis\'ement, si $2n+m>0$, on peut choisir des repr\'esentants pour param\'etrer notre
espace. Nous consid\'erons trois m\'ethodes~:
\vskip 0.5cm
{\bf Choix 1}~:

On choisit l'un des $z_i$ (disons $z_{j_0}=x_{j_0}+iy_{j_0}$) et on le place au point $i$ par une
transformation de $G$. Les autres points sont alors fix\'es~:

$$p_{j_0}=i,\quad p_j={z_j-x_{j_0}\over y_{j_0}},\quad q_l={t_l-x_{j_0}\over y_{j_0}}.$$
Dans ce cas, on param\`etre $C^+_{\{p_1,...,p_n\};\{q_1,...,q_m\}}$ par les coordonn\'ees des
$p_j=a_j+ib_j$ ($j\neq j_0$) et les $q_l$, l'orientation, dans ces coordonn\'ees de
$C^+_{\{p_1,...,p_n\};\{q_1,...,q_m\}}$ est celle donn\'ee par la forme~:

$$\Omega=\bigwedge_{j\neq j_0}(da_j\wedge db_j)\wedge dq_1\wedge...\wedge dq_m.$$
(L'ordre sur les indices $j$ n'importe pas car les 2-formes $da_j\wedge db_j$ commutent entre
elles).

\vskip 0.5cm
{\bf Choix 2}~:

On choisit l'un des $t_l$ (disons $t_{l_0}$) et on le place en 0 par une translation, puis on fait
une dilatation pour forcer le module de l'un des $z_j$ (disons $z_{j_0}$) \`a valoir 1~:

$$p_{j_0}={z_{j_0}-t_{l_0}\over \vert z_{j_0}-t_{j_0}\vert }=e^{i\theta_{j_0}},\quad p_j= {z_j-
t_{l_0}\over \vert z_{j_0}-t_{j_0}\vert},\quad q_{l_0}=0,\quad q_l={t_l-t_{l_0}\over \vert
z_{j_0}-t_{j_0}\vert }.$$

On param\`etre alors $C^+_{\{p_1,...,p_n\};\{q_1,...,q_m\}}$ par l'argument $\theta_{j_0}$
de $p_{j_0}$ (compris entre 0 et $\pi$) et par les coordonn\'ees des $p_j$ ($j\neq j_0$) et les
$q_l$ ($l\neq l_0$). L'orientation, dans ces coordonn\'ees de
$C^+_{\{p_1,...,p_n\};\{q_1,...,q_m\}}$ est celle donn\'ee par la forme~:

$$\Omega=(-1)^{l_0-1}d\theta_{j_0}\wedge\bigwedge_{j\neq j_0}(da_j\wedge db_j)\wedge
dq_1\wedge...\wedge \widehat{dq_{l_0}}\wedge...\wedge dq_m.$$
En effet, on part de la forme $\Omega$ du cas 1, avec $j_0=1$, puisque l'ordre des $p$
n'intervient pas, on place $q_{l_0}$ ``en t\^ete''~:

$$\Omega=(-1)^{l_0-1}dq_{l_0}\wedge da_2\wedge db_2\wedge...\wedge da_n\wedge
db_n\wedge dq_1\wedge...\wedge\widehat{dq_{l_0}}\wedge...\wedge dq_{l_1}\wedge...\wedge dq_n,$$
puis on effectue le changement de variables~:

$$p'_1={i-q_{l_0}\over\vert i-q_{l_0}\vert}=e^{i\theta_1},\quad p'_j={p_j-q_{l_0}\over\vert i-
q_{l_0}\vert} =a'_j+ib'_j \quad (2\leq j\leq n), \quad q'_k={q_k-q_{l_0}\over \vert i-
q_{l_0}\vert}\quad (k\neq l_0).$$
Dont le jacobien ${1\over (1+q^2_{l_0})^{m+n\over 2}}$ est strictement positif, pour obtenir la
forme annonc\'ee.

\vskip 0.5cm
{\bf Choix 3}~:

On choisit deux points $t_{l_0}<t_{l_1}$, on am\`ene par une translation le premier en $0$ et le
second en 1 par une dilatation.

$$p_j={z_j-t_{l_0}\over t_{l_1}-t_{l_0}},\quad q_{l_0}=0,\quad q_{l_1}=1,\quad q_l={t_l-
t_{l_0}\over t_{l_1}-t_{l_0}}.$$
On param\`etre $C^+_{\{p_1,...,p_n\};\{q_1,...,q_m\}}$ par les coordonn\'ees des
$p_j=a_j+ib_j$ et par les $q_l$ ($l\neq l_0$ et $l\neq l_1$). l'orientation est donn\'ee par la
forme~:

$$\Omega=(-1)^{l_0+l_1+1}\bigwedge_{j=1}^n(da_j\wedge db_j)\wedge
dq_1\wedge...\wedge
\widehat{dq_{l_0}}\wedge...\wedge\widehat{dq_{l_1}}\wedge...\wedge dq_m.$$
En effet, on part de la forme $\Omega$ du cas 1, on place $q_{l_0}$ et $q_{l_1}$ ``en t\^ete''~:

$$\Omega=(-1)^{l_0-1+l_1-2}dq_{l_0}\wedge dq_{l_1}\wedge da_2\wedge
db_2\wedge...\wedge da_n\wedge db_n\wedge
dq_1\wedge...\wedge\widehat{dq_{l_0}}\wedge...\wedge dq_n,$$
puis on effectue le changement de variables~:

$$p'_1={i-q_{l_0}\over q_{l_1}-q_{l_0}}=a_1+ib_1,\quad p'_j={p_j-q_{l_0}\over q_{l_1}-q_{l_0}}
=a'_j+ib'_j \quad (2\leq j\leq n), \quad q'k={q_k-q_{l_0}\over q_{l_1}-q_{l_0}}\quad (k\neq
l_0,~k\neq l_1).$$
Dont le jacobien ${q_{l_1}-q_{l_0}\over (q_{l_1}-q_{l_0})^{1+n+m}}$ est strictement positif, pour
obtenir la forme annonc\'ee.

\alinea{I.2. Compactification des espaces de configuration}

On plonge l'espace de configuration $C^+_{\{p_1,...,p_n\};\{q_1,...,q_m\}}$ dans une
vari\'et\'e compacte de la fa\c con suivante. Chaque fois que l'on prend deux points $A$ et $B$
du nuage de points $(z_j,\overline {z_j};t_l)$, on leur associe l'angle $Arg(B-A)$, \`a chaque triplet de points
$(A,B,C)$ du nuage, on associe l'\'el\'ement $[A-B,B-C,C-A]$ de l'espace projectif
$\BbbP^2(\BbbR)$ qu'ils d\'efinissent. On a ainsi une application~:

$$\tilde\Phi:Conf^+(z_j,t_l)\longrightarrow\BbbT^{(2n+m)(2n+m-1)}\times \left(\BbbP^2(\BbbR)
\right)^{(2n+m)(2n+m-1)(2n+m-2)}.$$
Cette application passe au quotient et il n'est pas difficile de montrer que l'on obtient ainsi un
plongement 

$$\Phi:C^+_{\{p_1,...,p_n\};\{q_1,...,q_m\}} \longrightarrow\BbbT^{(2n+m)(2n+m-1)}\times
\left(\BbbP^2(\BbbR) \right)^{(2n+m)(2n+m-1)(2n+m-2)}.$$
On d\'efinit la compactification $\overline{C^+_{\{p_1,...,p_n\};\{q_1,...,q_m\}}}$ de
$C^+_{\{p_1,...,p_n\};\{q_1,...,q_m\}}$ comme \'etant la fermeture dans $\BbbT^{(2n+m)(2n+m-
1)}\times \left(\BbbP^2(\BbbR) \right)^{(2n+m)(2n+m-1)(2n+m-2)}$ de
$\Phi\left(C^+_{\{p_1,...,p_n\};\{q_1,...,q_m\}}\right)$. On obtient ainsi une vari\'et\'e \`a coins
et on cherche son bord $\partial C^+_{\{p_1,...,p_n\};\{q_1,...,q_m\}}$.

Les points du bord s'obtiennent par une succession de collapses de points du nuage. On retrouve la description de M. Kontsevich \`a deux d\'etails pr\`es~: lorsque des
points a\'eriens (c'est \`a dire un ou des $p_j$) se rapprochent de $\BbbR$, il faut distinguer
entre quels $q_l$ ils arrivent, il y a trop de faces du bord, puisque les faces correspondant au
rapprochement de points terrestres (des $q_l$) non contigus est impossible sans que tous les
points qui les s\'eparent se rapprochent aussi. En codimension 1, on obtient deux types de faces~:
\alinea {I.2.1. Faces de type 1}
\qquad Parmi les points a\'eriens, $n_1$ points se rapprochent en un point $p$ qui reste a\'erien. Une
telle face existe si $n\geq n_1\geq 2$. A la limite, on obtient une vari\'et\'e produit~:

$$F=\partial_{\{p_{i_1},...,p_{i_{n_1}}\}} C^+_{\{p_1,...,p_n\};\{q_1,...,q_m\}} =
C_{\{p_{i_1},...,i_{n_1}\}}  \times C_{\{p,p_1,...,\widehat{p_{i_1}},...,\widehat{p_{i_{n_1}}},
...,p_n\};\{q_1,...,q_m\}}\eqno{(*)}$$
o\`u l'espace $C_{\{p_1,...,p_{n_1}\}} $ est le quotient de l'espace $Conf(z_1,..., z_{n_1})$ par
l'action du groupe $G'$ des transformations $z_j\mapsto az_j+b$ ($a>0$ et $b\in\BbbC$). c'est
une vari\'t\'e de dimension $2n_1-3$ ($n_1\geq 2$). On la plonge dans un produit de tores et
d'espaces projectifs comme pour $C^+_{\{p_1,...,p_n\};\{q_1,...,q_m\}}$. Enfin on l'oriente de la
fa\c con suivante; $z_1$ est plac\'e en 0 par une translation complexe puis $\vert z_2\vert$ est
normalis\'e \`a 1 par une dilatation,

$$p_1=0,\quad p_2={z_2-z_1\over\vert z_2-z_1\vert}=e^{i\theta_2},\quad p_j={z_j-z_1\over
\vert z_2-z_1\vert}=a_j+ib_j$$
et on prend l'orientation d\'efinie par la forme~:

$$\Omega_1=d\theta_2\wedge\bigwedge_{j\geq 3}(da_j\wedge db_j).$$

Orientons maintenant la face $F$. On choisit la forme volume $\Omega_1\wedge\Omega_2$
sur le produit $(*)$ o\`u $\Omega_2$ est l'une des formes d\'efinies ci-dessus pour orienter
$C^+_{\{p_j\};\{q_l\}}$. L'orientation de la face \`a partir de celle de $\Omega$ est $\pm
\Omega_1\wedge\Omega_2$.
\lemme{I.2.1}
La face $F$ est orient\'ee par $\Omega_F=-\Omega_1\wedge\Omega_2$.
\dem
On a vu que l'on pouvait changer l'ordre des points $p_j$ de $C^+_{\{p_j\};\{q_l\}}$ sans
changer l'orientation. On renum\'erote les points $p_{i_1}$,..., $p_{i_{n_1}}$ en $p_1$,
$p_2$,..., $p_{n_1}$, puis on fixe $p_1=i$~:

$$\Omega=\bigwedge_{j=2}^{n_1}(da_j\wedge db_j)\wedge \Omega_2,$$
ensuite on change de variables dans le premier facteur en posant~:

$$p'_2=e^{i\theta_2},\quad p'_j={p_j\over \vert p_2-i\vert}=a'_j+ib'_j\qquad (j=3,...,n_1).$$
Lorsque les $n_1$ premiers points collapsent, on agrandit le petit nuage qu'ils forment en
normalisant la distance qui s\'epare les 2 premiers \`a 1. Posons $\rho_2=\vert p_2-i\vert$. le
changement de variable donne pour $\Omega$ la forme~:

$$\Omega'=d\rho_2\wedge d\theta_2\wedge\bigwedge_{j\geq 3}(da'_j\wedge db'_j)\wedge
\Omega_2.$$
La face est obtenue lorsque $\rho_2\rightarrow 0$. Or $\rho_2>0$, on doit donc l'orienter avec

$$\Omega_F=-d\theta_2\wedge\bigwedge_{j\geq 3}(da'_j\wedge db'_j)\wedge \Omega_2=-
\Omega_1 \wedge \Omega_2.$$
\qed
\alinea{I.2.2. Face de type 2}
\qquad Parmi les points du nuage, $n_1$ points a\'eriens et $m_1$ points terrestres se rapprochent en
un point $q$ terrestre. Une telle face existe si $n+m>n_1+m_1$ et $2n_1+m_1\geq 2$. A la limite, on obtient une vari\'et\'e produit~:

$$\eqalign{F&=\partial_{\{p_{i_1},...,p_{i_{n_1}}\};\{q_{l+1},...,q_{l+m_1}\}} C^+_{\{p_1,...,
p_n\}; \{q_1,...,q_m\}} \cr &= C_{\{p_{i_1},...,p_{i_{n_1}}\};\{ q_{l+1},...,q_{l+m_1}\}}\times
C_{\{p_1,...,\widehat{p_{i_1}}, ...,\widehat{p_{i_{n_1}}}, ...,p_n\};
\{q_1,...,q_l,q,q_{l+m_1+1},...,q_m\}}.}\eqno{(*)}$$
On appelle $\Omega_1$ et $\Omega_2$ l'une des formes volumes de chacun des facteurs de
ce produit. La forme $\Omega_1\wedge \Omega_2$ est une forme volume sur $F$. On donne
l'orientation de $F$ \`a partir de celle de l'espace de configuration de d\'epart en terme de cette
forme.
\lemme{I.2.2}
Avec nos notations, la face $F$ est orient\'ee par~:
$$\Omega_F=(-1)^{lm_1+l+m_1}\Omega_1\wedge\Omega_2.$$
\dem
Il faut consid\'erer six types de nuages diff\'erents~:
\vskip 0.5cm
{\bf Sous-cas 1}~: $n>n_1>0$
\dessin{25mm}{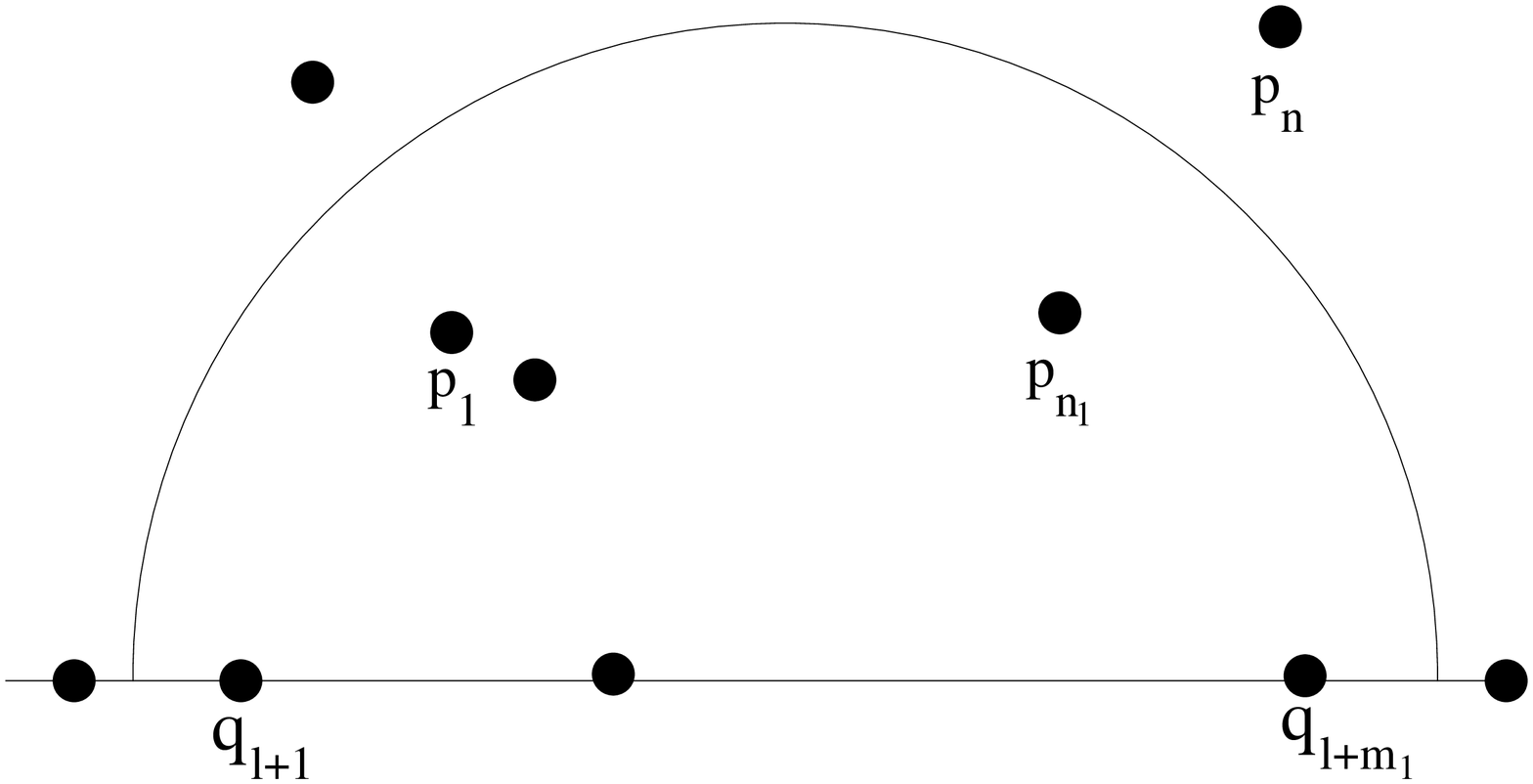}
On suppose que $p_1$,..., $p_{n_1}$ et $q_{l+1}$, ..., $q_{l+m_1}$ collapsent. On param\`etre
l'espace $C^+_{\{p_1,..., p_n\}; \{q_1,...,q_m\}}$ par $p_{n_1+1}=i$. La forme d'orientation est~:

$$\eqalign{\Omega=&(-1)^{lm_1}da_1\wedge db_1\wedge...\wedge da_{n_1}\wedge
db_{n_1}\wedge dq_{l+1}\wedge...\wedge dq_{l+m_1}\bigwedge \cr
&\hskip 3cm\bigwedge da_{n_1+2}\wedge db_{n_1+2}\wedge...\wedge da_n\wedge
db_n\wedge dq_1\wedge...\wedge dq_l\wedge dq_{l+m_1+1}\wedge...\wedge dq_m.\cr
=&(-1)^{lm_1+l+m_1+1}db_1\wedge...\wedge da_{n_1}\wedge db_{n_1}\wedge
dq_{l+1}\wedge...\wedge dq_{l+m_1}\bigwedge \cr
&\hskip 3cm\bigwedge da_{n_1+2}\wedge db_{n_1+2}\wedge...\wedge da_n\wedge
db_n\wedge dq_1\wedge...\wedge dq_l\wedge da_1\wedge dq_{l+m_1+1}\wedge...\wedge
dq_m. }$$
(Certains termes peuvent ne pas appara\^\i tre, par exemple si $n=n_1+1$ ou $m=m_1$). On
change de variables en posant~:

$$a_1=q,\quad p'_j={p_j-a_1\over b_1}\quad (2\leq j\leq n_1),\qquad q'_k={q_k-a_1\over
b_1} \quad (l+1\leq k\leq l+m_1).$$
Alors on peut \'ecrire de fa\c con un peu abusive~:

$$\Omega=(-1)^{lm_1+l+m_1+1}db_1\bigwedge \Omega_1\wedge\Omega_2$$
et puisque $b_1>0$ et la face $F$ est obtenue pour $b_1=0$, son orientation est donn\'ee par~:

$$\Omega_F=(-1)^{lm_1+l+m_1}\Omega_1\wedge\Omega_2.$$

\vskip 0.5cm
{\bf Sous-cas 2}~: $n=n_1>0$ (et donc $m\ge m_1+1$) et $l>0$
\dessin{25mm}{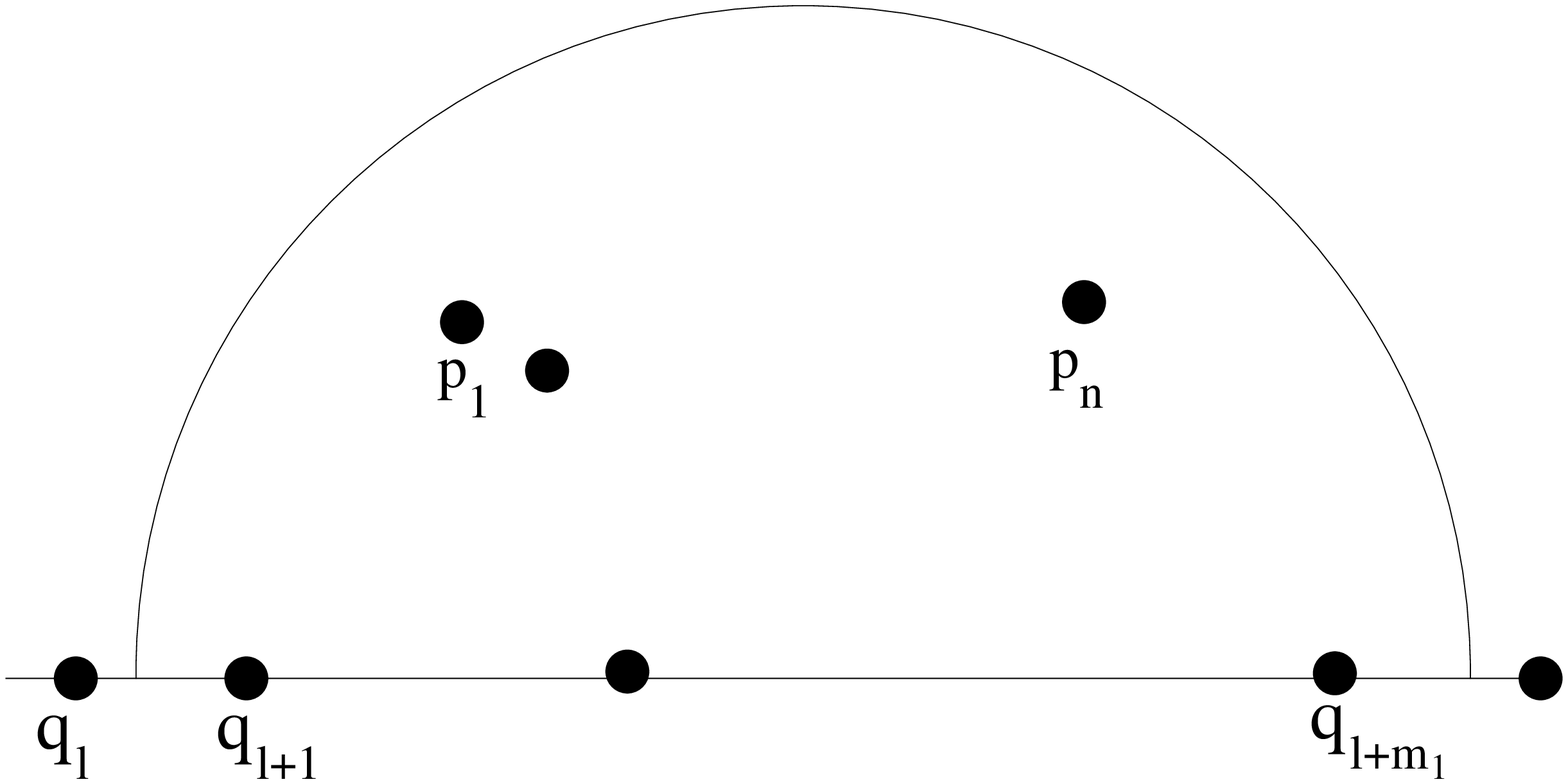}
On suppose que $p_1$,..., $p_n$ et $q_{l+1}$, ..., $q_{l+m_1}$ collapsent. On param\`etre
l'espace $C^+_{\{p_1,..., p_n\}; \{q_1,...,q_m\}}$ par $q_l=0$, $q_{l+1}=1$. La forme
d'orientation est~:

$$\eqalign{\Omega&=(-1)^{2l+2}da_1\wedge db_1\wedge...\wedge da_n\wedge db_n\wedge
dq_1\wedge...\wedge \widehat{dq_l}\wedge \widehat{dq_{l+1}}\wedge...\wedge dq_m\cr
&=(-1)^{(l-1)(m_1-1)} da_1\wedge db_1\wedge...\wedge da_n\wedge db_n\wedge
dq_{l+2}\wedge ...\wedge dq_{l+m_1}\bigwedge\cr
&\hskip 4cm\bigwedge dq_1\wedge...\wedge dq_{l-1}\wedge dq_{l+m_1+1}\wedge...\wedge
dq_m\cr}.$$
On change de variables en posant~:

$$a_1+ib_1-1=\rho_1e^{i\theta_1}\quad p'_j={p_j-1\over \rho_1}\quad (2\leq j\leq n),\qquad
q'_k={q_k-1\over \rho_1} \quad (l+2\leq k\leq l+m_1).$$
Alors

$$\eqalign{\Omega&= (-1)^{lm_1+l+m_1+1}\rho_1^{-(n+m_1-5)}d\rho_1\wedge d\theta_1\wedge
da_2\wedge db_2\wedge...\wedge da_n\wedge db_n\wedge dq_{l+2}\wedge ...\wedge
dq_{l+m_1}\bigwedge\cr
&\hskip 4cm\bigwedge dq_1\wedge...\wedge dq_{l-1}\wedge dq_{l+m_1+1}\wedge...\wedge
dq_m\cr}.$$

On peut donc \'ecrire de fa\c con un peu abusive~:

$$\Omega\simeq (-1)^{lm_1+l+m_1+1}d\rho_1\bigwedge \Omega_1\wedge(-1)^{l-1+l+1-2}
\Omega_2$$
et puisque $\rho_1>0$ et la face $F$ est obtenue pour $\rho_1=0$, son orientation est
donn\'ee par~:

$$\Omega_F=(-1)^{lm_1+l+m_1}\Omega_1\wedge\Omega_2.$$

\vskip 0.5cm
{\bf Sous-cas 3}~: $n=n_1>0$ (et donc $m\ge m_1+1$) et $l=0$
\dessin{25mm}{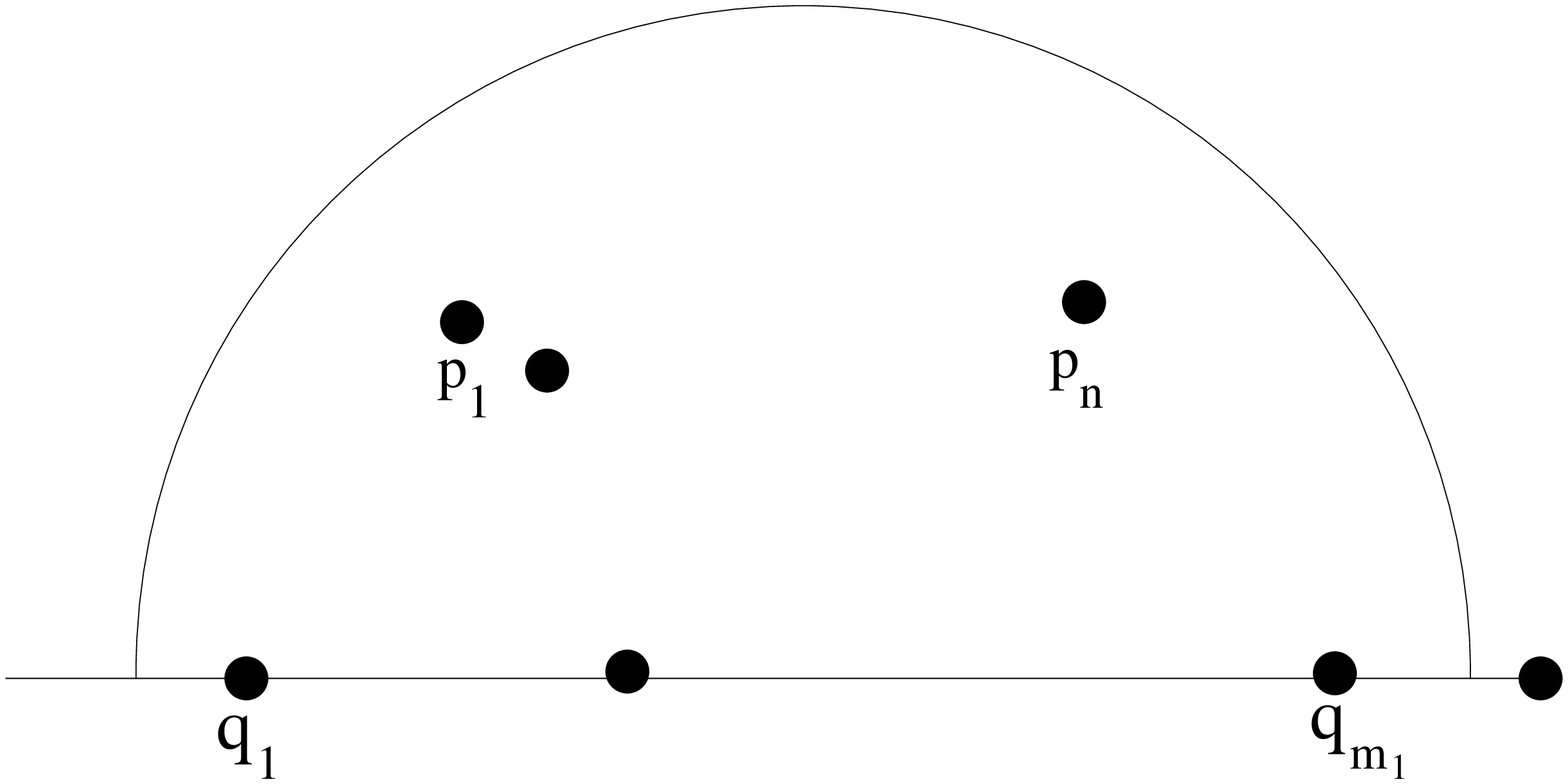}
C'est le m\^eme calcul que ci-dessus, on pose $q_{m_1}=0$, $q_{m_1+1}=1$, on obtient~:

$$\Omega=(-1)^{2m_1+2}da_1\wedge db_1\wedge...\wedge da_n\wedge db_n\wedge
dq_1\wedge...\wedge \widehat{dq_{m_1}}\wedge \widehat{dq_{m_1+1}}\wedge...\wedge
dq_m.$$
On change de variables en posant~:

$$a_1+ib_1=\rho_1e^{i\theta_1}\quad p'_j={p_j\over \rho_1}\quad (2\leq j\leq n),\qquad
q'_k={q_k\over \rho_1} \quad (1\leq k\leq m_1-1).$$
Alors

$$\eqalign{\Omega&=\rho_1d\rho_1\wedge d\theta_1\wedge da_2\wedge
db_2\wedge...\wedge da_n\wedge db_n\wedge dq_1\wedge ...\wedge dq_{m_1-1}
\bigwedge dq_{m_1+2}\wedge...\wedge dq_m\cr 
&\simeq d\rho_1\bigwedge (-1)^{m_1-1}\Omega_1\wedge (-1)^{1-1+2-2}
\Omega_2 \cr}$$
et puisque $\rho_1>0$ et la face $F$ est obtenue pour $\rho_1=0$, son orientation est
donn\'ee par~:

$$\Omega_F=(-1)^{lm_1+l+m_1}\Omega_1\wedge\Omega_2.$$

\vskip 0.5cm
{\bf Sous-cas 4}~: $n\ge n_1=0$ (et donc $m_1>1$) et $l>0$
\dessin{25mm}{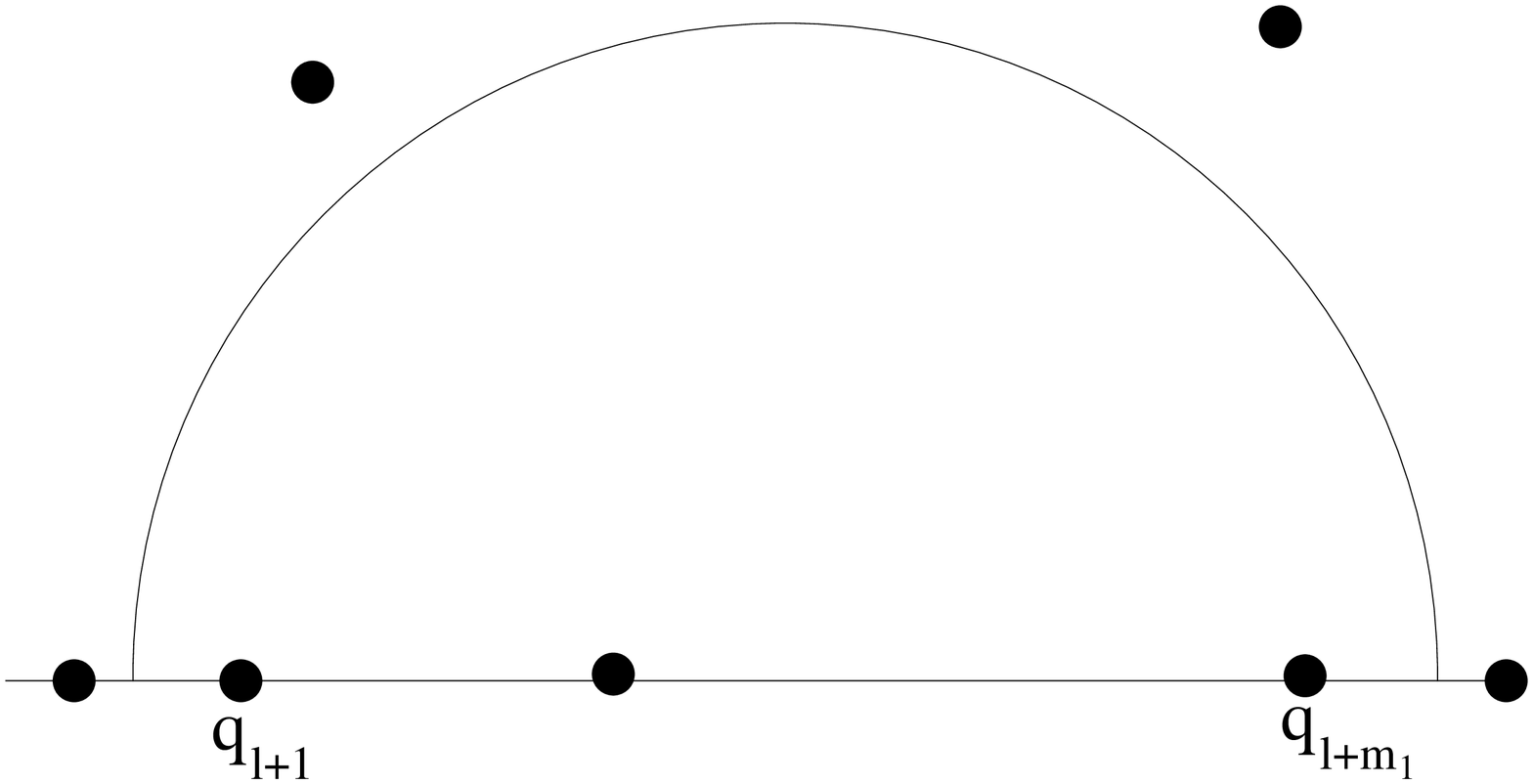}
On suppose que $q_{l+1}$, ..., $q_{l+m_1}$ collapsent. On param\`etre l'espace
$C^+_{\{p_1,..., p_n\}; \{q_1,...,q_m\}}$ par $q_l=0$, $q_{l+1}=1$. La forme d'orientation est~:

$$\eqalign{\Omega&=(-1)^{2l+2}da_1\wedge db_1\wedge...\wedge da_n\wedge db_n\wedge
dq_1\wedge...\wedge \widehat{dq_l}\wedge \widehat{dq_{l+1}}\wedge...\wedge dq_m\cr
&=(-1)^{(l-1)(m_1-1)} dq_{l+2}\wedge ...\wedge dq_{l+m_1}\bigwedge\cr
&\hskip 4cm\bigwedge da_1\wedge db_1\wedge...\wedge da_n\wedge db_n\wedge
dq_1\wedge...\wedge dq_{l-1}\wedge dq_{l+m_1+1}\wedge...\wedge dq_m\cr}.$$
On change de variables en posant~:

$$q'_k={q_k-q_{l+2}\over q_{l+2}-1} \quad (l+3\leq k\leq l+m_1).$$
Alors~:

$$\eqalign{\Omega&\simeq (-1)^{lm_1+l+m_1+1}dq_{l+2}\wedge  dq'_{l+3}\wedge ...\wedge
dq'_{l+m_1}\bigwedge\cr
&\hskip 4cm\bigwedge da_1\wedge db_1\wedge...\wedge da_n\wedge db_n\wedge
dq_1\wedge...\wedge dq_{l-1}\wedge dq_{l+m_1+1}\wedge...\wedge dq_m\cr}.$$

On peut donc \'ecrire de fa\c con un peu abusive~:

$$\Omega=(-1)^{lm_1+l+m_1+1}dq_{l+2}\bigwedge \Omega_1\wedge(-1)^{l-1+l+1-2}
\Omega_2$$
et puisque $q_{l+2}-1>0$ et la face $F$ est obtenue pour $q_{l+2}-1=0$, son orientation est
donn\'ee par~:

$$\Omega_F=(-1)^{lm_1+l+m_1}\Omega_1\wedge\Omega_2.$$

\vskip 0.5cm
{\bf Sous-cas 5}~: $n\ge n_1=0$ (et donc $m_1>1$), $l=0$ et $m_1<m$
\dessin{25mm}{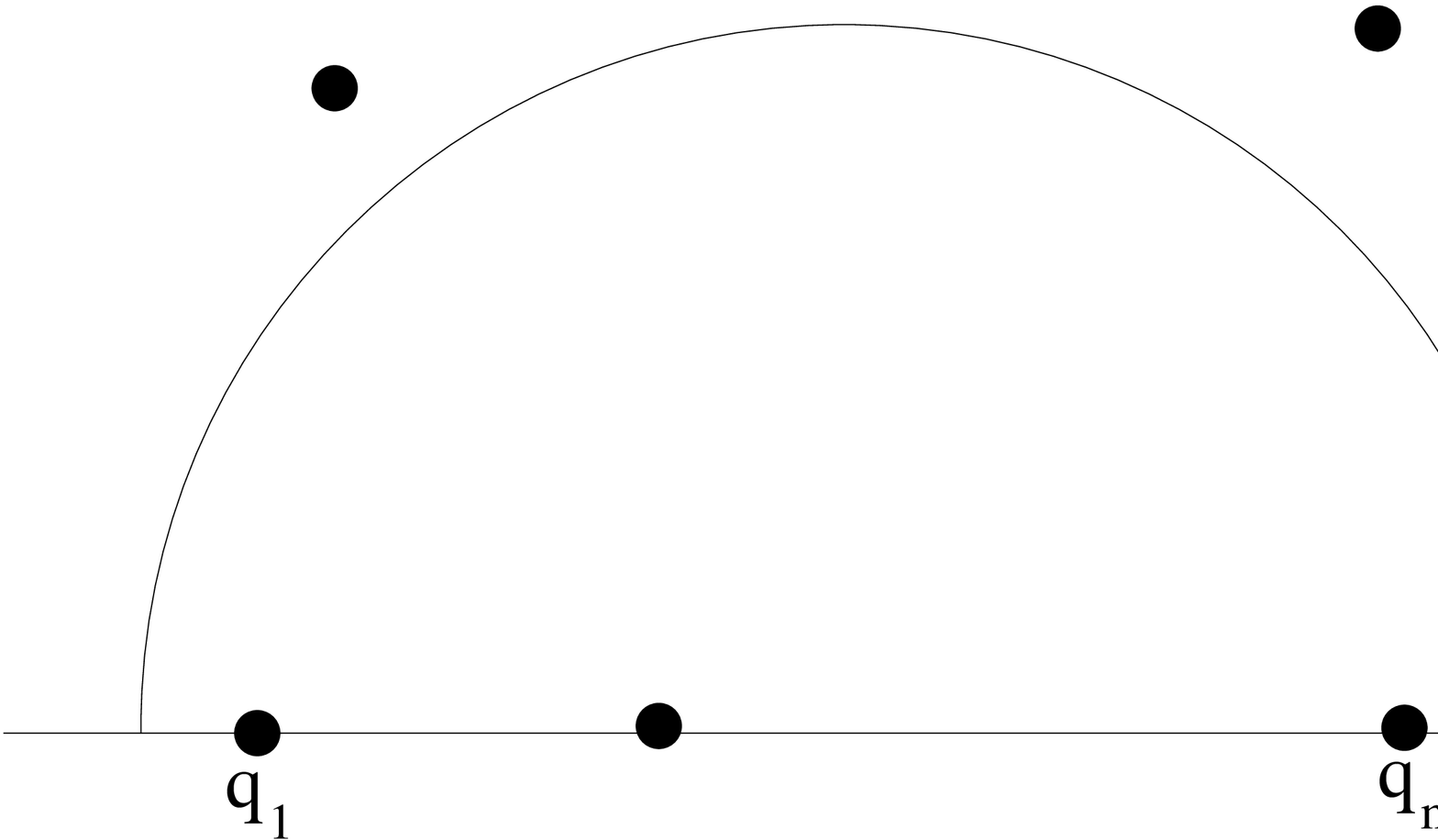}
On pose $q_{m_1}=0$, $q_{m_1+1}=1$, on obtient~:

$$\Omega=(-1)^{2m_1+2}da_1\wedge db_1\wedge...\wedge da_n\wedge db_n\wedge
dq_1\wedge...\wedge \widehat{dq_{m_1}}\wedge \widehat{dq_{m_1+1}}\wedge...\wedge
dq_m.$$
On change de variables en posant~:

$$q'_k={q_k-q_{m_1-1}\over -q_{m_1-1}}\quad (1\leq k\leq m_1-2)$$
Alors~:

$$\eqalign{\Omega&\simeq dq'_1\wedge...\wedge dq'_{m_1-2}\bigwedge dq_{m_1-1}
\bigwedge da_1\wedge db_1\wedge...\wedge da_n\wedge db_n\wedge dq_{m_1+2} \wedge
...\wedge dq_m\cr
&=(-1)^{m_1}dq_{m_1-1} \wedge \Omega_1\wedge (-1)^{1+2-1}\Omega_2.}$$
Maintenant $q_{m_1-1}<0$ et l'orientation de la face est encore~:

$$\Omega_F=(-1)^{lm_1+m_1+l}\Omega_1\wedge \Omega_2.$$

\vskip 0.5cm
{\bf Sous-cas 6}~: $n\ge n_1=0$ (et donc $m_1>1$), $l=0$ et $m_1=m$ et donc $n>0$
\dessin{25mm}{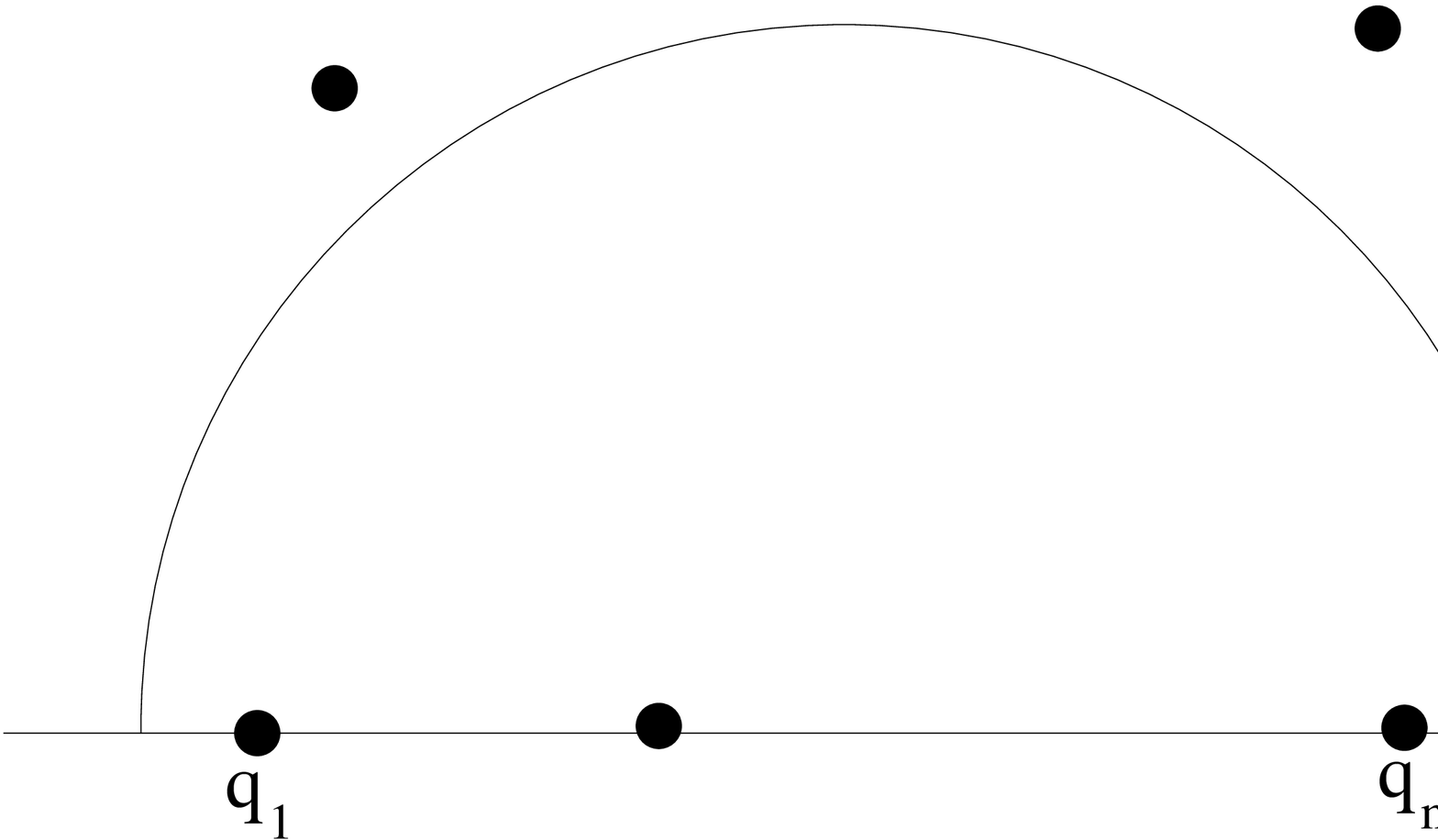}
On pose $q_1=0$ et $p_1=e^{i\theta_1}$. La forme $\Omega$ est

$$\eqalign{\Omega&=d\theta_1\wedge da_2\wedge db_2\wedge...\wedge da_n\wedge
db_n\wedge dq_2 \wedge...\wedge dq_m\cr
&=(-1)^{m-1}dq_2\wedge...\wedge dq_m\bigwedge d\theta_1\wedge da_2\wedge
db_2\wedge...\wedge da_n\wedge db_n.}$$
On change de variables en posant~:

$$q'_k={q_k-q_2\over q_2}\quad (3\leq k\leq m)$$
Alors~:

$$\Omega\simeq (-1)^{m-1} dq_2\wedge (-1)^{1-1+2-2}\Omega_1\wedge \Omega_2.$$
Puisque $q_2>0$ et $F$ appara\^\i t pour $q_2=0$, l'orientation de la face est encore~:

$$\Omega_F=(-1)^{lm_1+l+m_1}\Omega_1\wedge \Omega_2.$$
\vskip 0.5cm
Tout nuage de point correspondant \` a une face de type 2 rel\`eve d'un de ces six sous-cas. Ceci termine la d\'emonstration du lemme I.2.2.
\qed
\paragraphe{II. Alg\`ebres sym\'etriques et ext\'erieures sur les espaces gradu\'es}
\alinea{II.1. La cat\'egorie des espaces gradu\'es}
Un espace vectoriel sur un corps $k$ est {\sl gradu\'e\/} s'il est muni d'une $\Z-$graduation~:
$$V=\bigoplus_{n\in\Z}V_n$$
Le degr\'e d'un \'el\'ement homog\`ene $x$ sera not\'e $|x|$. Un espace gradu\'e sera toujours consid\'er\'e comme un super-espace vectoriel, la $\Z_2-$graduation \'etant d\'eduite de la $\Z-$graduation~:
$$V_+=\bigoplus_{n\in\Z}V_{2n}\hbox to 12mm{} V_-=\bigoplus_{n\in\Z}V_{2n+1}$$
Si $V$ et $W$ sont des espaces gradu\'es, il existe une graduation naturelle sur $V\oplus W$, $V\otimes W$, $\mop{Hom}_k(V,W)$. Un {\sl morphisme d'espaces gradu\'es\/} entre $V$ et $W$ est par d\'efinition un \'el\'ement de {\sl degr\'e z\'ero\/} dans $\mop{Hom}_k(V,W)$.
\ssq
Une {\sl alg\`ebre gradu\'ee\/} est un espace gradu\'e $B$ muni d'une structure d'alg\`ebre telle que la multiplication $m:B\otimes B\rightarrow B$ est un morphisme d'espaces gradu\'es, c'est-\`a-dire~:
$$B_iB_j\subset B_{i+j}$$
On d\'efinit de la m\^eme mani\`ere les notions de $B-$ modules gradu\'es \`a gauche ou \`a droite. Si $A$ et $B$ sont deux alg\`ebres gradu\'ees, le produit~:
$$\eqalign{m_{A\otimes B}:A\otimes B\otimes A\otimes B  &\longrightarrow
                                                A\otimes B\cr
        a\otimes b\otimes a'\otimes b'                  &\longmapsto
                                                (-1)^{|b||a'|}aa'\otimes bb'\cr}$$
est associatif et munit $A\otimes B$ d'une structure d'alg\`ebre gradu\'ee. Si $M$ (resp. $N$) est un $A$-module (resp. un $B$-module) gradu\'e \`a gauche, la m\^eme r\`egle des signes (la {\sl r\`egle de Koszul\/}) permet de d\'efinir une structure de $A\otimes B-$module gradu\'e \`a gauche sur $M\otimes N$.
\ssq
Si $A,A',B,B'$ sont des espaces gradu\'es, l'identification de $\mop{Hom}_k(A\otimes B, A'\otimes B')$ avec $\mop{Hom}_k(A,A')\otimes\mop{Hom}_k(B,B')$ se fait avec la m\^eme r\`egle des signes~:
$$(f\otimes g)(a\otimes b)=(-1)^{|g||a|}f(a)\otimes g(b)$$
\ssq
Une alg\`ebre gradu\'ee est dite {\sl commutative\/} si on a~:
$$xy-(-1)^{|x||y|}yx=0$$
Une cog\`ebre gradu\'ee $C$ se d\'efinit de mani\`ere similaire~: la comutiplication doit v\'erifier~:
$$\Delta C_j\subset \sum_{k+l=j}C_k\otimes C_l$$
On d\'efinit une structure de cog\`ebre gradu\'ee sur le produit tensoriel de deux cog\`ebres gradu\'ees en appliquant la m\^eme r\`egle sur les signes que dans le cas des alg\`ebres.
\ssq
Une {\sl d\'erivation de degr\'e $i$} dans une alg\`ebre gradu\'ee $B$ est un morphisme lin\'eaire $d:B\rightarrow B$ de degr\'e $i$ tel que~:
$$d(xy)=dx.y+(-1)^{i|x|}x.dy$$
ce qui s'\'ecrit encore~:
$$dm=m(d\otimes I+I\otimes d)$$
o\`u $m$ d\'esigne la multiplication de l'alg\`ebre (attention \`a la r\`egle des signes).
Une {\sl cod\'eriva\-tion de degr\'e $i$} dans une cog\`ebre gradu\'ee $C$ est un morphisme lin\'eaire $d:C\rightarrow C$ de degr\'e $i$ tel que si $\Delta X=\sum_{(x)}x'\otimes x''$ on a~:
$$\Delta dx=\sum_{(x)}dx'\otimes x''+(-1)^{i|x'|}x'\otimes dx''$$
ou encore~:
$$\Delta d=(d\otimes I + I\otimes d)\Delta$$
Enfin une {\sl alg\`ebre de Lie gradu\'ee\/} est un espace vectoriel gradu\'e $\g g$ muni d'un crochet $[.,.]$ tel que~:
\ssq
1) $[\g g_i,\g g_j]\subset \g g_{i+j}$
\ssq
2) $[x,y]=-(-1)^{|x||y|}[y,x]$
\def\toto #1#2#3{(-1)^{|#1||#3|}[[#1,#2],#3]}
\ssq
3) $\permuc\toto xyz =0$

\qquad(identit\'e de Jacobi gradu\'ee).
\smallskip
L'identit\'e de Jacobi gradu\'ee s'exprime aussi en disant que $\mop{ad}x=[x,.]$ est une d\'erivation (de degr\'e $|x|$).
\ssq
Une alg\`ebre de lie gradu\'ee est {\sl diff\'erentielle} si elle est munie d'une diff\'erentielle
$d$ de degr\'e 1 ($d:{\frak g}\longrightarrow{\frak g}[1]$), telle que~:

$$d^2=0,\quad d\left([x,y]\right)=\left[ dx,y\right]+(-1)^{1.\vert x\vert}\left[x,dy\right].$$
\alinea{II.2 La r\`egle de Koszul}
La raison profonde qui fait que ``la r\`egle des signes marche'' est la suivante~: la cat\'egorie des espaces vectoriels $\Z_2-$gradu\'es munie du produit tensoriel $\otimes$ usuel et des applications~:
$$\eqalign{\tau_{A,B}:A\otimes B        &\longrightarrow B\otimes A\cr
                a\otimes b              &\longmapsto (-1)^{|a||b|}b\otimes a\cr}$$
est une {\sl cat\'egorie tensorielle tress\'ee\/}, c'est \`a dire que les {\sl tressages\/} $\tau_{A,B}$ sont fonctoriels~:
$$\matrix{A\otimes B    &\fleche {12}_{\tau_{A,B}}      &B\otimes A\cr
        ^{f\otimes g}\big\downarrow \hbox to 5mm{}&                     &\hbox to 5mm{}\big\downarrow ^{g\otimes f}\cr
        A'\otimes B'    &\fleche {12}_{\tau_{A',B'}}    &B'\otimes A'\cr}$$
et v\'erifient~:
$$\tau_{A\otimes B,C}=(\tau_{A,C}\otimes I_B)(I_A\otimes \tau_{B,C})$$
De ces deux propri\'et\'es on d\'eduit facilement l'\'equation de l'hexagone, c'est-\`a-dire la commutativit\'e du diagramme suivant~:
$$\matrix{&     A\otimes B\otimes C &\cr
        \ \swarrow      &&\!\!\searrow  \cr
        A\otimes C\otimes B &&  B\otimes A\otimes C\cr
        \big\downarrow          &&      \big\downarrow\cr
        C\otimes A\otimes B     &&      B\otimes C\otimes A\cr
        \ \searrow      &&\!\!\swarrow  \cr
        &       C\otimes B\otimes A     &\cr}$$
La cat\'egorie tensorielle tress\'ee des espaces $\Z_2-$gradu\'es peut aussi se voir comme la cat\'ego\-rie des modules sur l'alg\`ebre de Hopf quasi-triangulaire $(H_2, R)$ o\`u $H_2$ est l'alg\`ebre du groupe $\Z_2$ munie de la multiplication et de la comultiplication usuelle, mais o\`u la $R-$matrice est non triviale.
\ssq
Dans cette cat\'egorie le carr\'e des tressages est toujours l'identit\'e (c'est une cat\'egorie tensorielle stricte). On peut faire de m\^eme avec des espaces $\Z_k-$gradu\'es en rempla\c cant $-1$ par $e^{{2i\pi \over k}}$. On obtient ainsi la cat\'egorie des {\sl espaces vectoriels anyoniques\/}, qui est tress\'ee de mani\`ere effective pour $k\ge 3$ [M].
\alinea{II.3. D\'ecalages}
Soit $V$ un espace gradu\'e. On pose~:
$$V[1]=V\otimes k[1]$$
o\`u $k[1]$ est l'espace gradu\'e tel que $k_n=\{0\}$ pour $n\not =-1$ et $k_{-1}=k$. Autrement dit $V[1]$ et $V$ ont m\^eme espace vectoriel sous-jacent, mais le degr\'e d'un \'el\'ement est baiss\'e d'une unit\'e dans $V[1]$. On posera en outre~:
$$[n]=[1]^n$$
pour tout entier $n$.
\alinea{II.4. Alg\`ebres sym\'etriques et ext\'erieures}
L'alg\`ebre sym\'etrique $S(V)$ (resp. l'alg\`ebre ext\'erieure $\Lambda (V)$) est d\'efinie par~:
$$S(V)=T(V)/<x\otimes y-(-1)^{|x||y|}y\otimes x> \hbox to 15mm{\hfill resp.}
\Lambda(V)=T(V)/<x\otimes y+(-1)^{|x||y|}y\otimes x>$$
Ce sont des espaces gradu\'es de mani\`ere naturelle. La proposition suivante est implicite dans \cite {K1}~:
\prop{II.4.1 \rm (sym\'etrisation)}
Pour tout espace vectoriel gradu\'e $V$ et pour tout $n>0$ on a un isomorphisme naturel~:
$$\Phi_n:S^n(V[1])\widetilde{\fleche 8} \Lambda^n(V)[n]$$
donn\'e par~:
$$\Phi_n(x_1.\ldots .x_n)=\alpha(\uple xn)x_1\wedge\cdots\wedge x_n$$
o\`u, pour des $x_i$ homog\`enes, $\alpha(\uple xn)$ d\'esigne la signature de la permutation ``unshuffle'' qui range les $x_i$ pairs dans $V$ \`a gauche sans les permuter, et les $x_i$ impairs dans $V$ \`a droite sans les permuter.
\dem
Soit $I$ (resp. $J$) l'ensemble des $i$ tels que $x_i$ soit de degr\'e pair (resp. impair), et $\alpha(I,J)=\alpha(\uple xn)$ la signature de la permutation-rangement associ\'ee. L'isomorphisme $\Phi_n$ est donn\'e par la restriction \`a $S^n(V[1])$ de la composition des trois fl\`eches du diagramme ci-dessous (la fl\`eche sup\'erieure est un isomorphisme d'alg\`ebres)~:
$$\matrix
{S(V[1])        &\widetilde{\fleche{12}}        &S(V[1]_+)\otimes S(V[1]_-)\cr
x_1\ldots x_n   &\longmapsto            &x_I\otimes x_J \cr
&&                                      \big\downarrow  \cr
\Lambda(V)      &\widetilde{\fleche{12}}        &\Lambda(V_-)\otimes \Lambda(V_+)\cr
\alpha(I,J)x_1\wedge\cdots\wedge x_n       &\longmapsto
&x_{\wedge I}\otimes x_{\wedge J}\cr}$$
Enfin si les $x_j$ sont de degr\'e $d_j$ dans $V[1]$, $x_1\ldots x_n$ est de degr\'e $d_1+\cdots +d_n$ dans $S^n(V[1])$, donc de degr\'e $d_1+\cdots +d_n+n$ dans $S^n(V)$. $\Phi_n(x_1\ldots x_n)$ est donc de degr\'e $d_1+\cdots +d_n+n$ dans $\Lambda^n(V)$, donc de degr\'e  $d_1+\cdots +d_n$ dans $\Lambda^n(V)[n].$
\qed
{\sl Remarque\/}~: L'application $\Phi=\oplus\Phi_n$ est un morphisme d'espace
vectoriel gradu\'e mais pas d'alg\`ebre. Il est d'ailleurs vain de vouloir chercher un isomorphisme d'alg\`ebres entre $S(V[1])$ et $\bigoplus \Lambda^n(V)[n]$, car deux \'el\'ements de parit\'e oppos\'ee commutent dans le premier cas, et anticommutent dans le second cas.
\alinea{II.5. Un exemple : $Tens(\R^d)$}
\qquad L'alg\`ebre des tenseurs contravariants totalement antisym\'etriques est une alg\`ebre
naturellement gradu\'ee par l'ordre des tenseurs. On aimerait la voir comme l'espace sous-jacent
\`a une alg\`ebre sym\'etrique. Notons donc $V$ l'espace vectoriel ${\cal X}
\left(\BbbR^d\right)$ des champs de vecteurs sur $\BbbR^d$, gradu\'e par $V=V_0$. On
identifie $Tens_n\left(\BbbR^d\right)=\wedge^nV$ \`a $S^n\left(V[1]\right)[-n]$ par $\Phi_n$.

Dans la suite, on posera

$$T_{poly}\left(\BbbR^d\right)=Tens\left(\BbbR^d\right)[1].$$

\paragraphe{III. Vari\'et\'es formelles gradu\'ees}
\alinea{III.1. Vari\'et\'es formelles}
On se place sur le corps des r\'eels ou des complexes. On se donne un voisinage ouvert $U$ de $0$ dans $\R^d$. Une fonction analytique $\varphi$ sur $U$ \`a valeurs dans $\C$ est d\'etermin\'ee par son d\'eveloppement de Taylor en $0$~:
$$\varphi(x)=\sum_{\alpha\in\N^d}{x^\alpha\over{\alpha !}}(\partial^\alpha\varphi)(0)$$
On \'etablit ainsi une dualit\'e non d\'eg\'en\'er\'ee entre les fonctions analytiques sur $U$ et les distributions de support $\{0\}$. Plus abstraitement on peut remplacer les fonctions analytiques par les jets d'ordre infini au point $0$.
\ssq
On appelle {\sl vari\'et\'e formelle\/}, ou voisinage formel de $0$, l'espace $\overline\Cal C$ des distributions de support $\{0\}$. La structure d'alg\`ebre commutative sur l'espace des fonctions analytiques sur $U$ d\'etermine une structure de cog\`ebre cocommutative sur son dual restreint, qui est exactement $\overline\Cal C$. La comultiplication est donn\'ee par~:
$$<\Delta v,\varphi\otimes\psi>=<v,\varphi\psi>.$$
Consid\'erant l'espace tangent $V$ \`a la vari\'et\'e $U$ en $0$, on a en fait un isomorphisme de cog\`ebres entre $\overline\Cal C$ et $S(V)$, o\`u la comultiplication $\Delta$ de $S(V)$ est le morphisme d'alg\`ebres tel que $\Delta(v)=v\otimes 1+1\otimes v$ pour $v\in V$.
\ssq
On consid\`erera la version point\'ee~:
$$\Cal C=S^+(V)=\bigoplus_{n\ge 1}S^n(V)$$
C'est la cog\`ebre colibre cocommutative sans co-unit\'e construite sur $V$. C'est aussi le dual restreint de l'alg\`ebre des jets d'ordre infini qui s'annulent en $0$. On remarque que $\Delta v=0$ si et seulement si $v$ appartient \`a $V$.
\ssq
Un {\sl champ de vecteurs\/} sur la vari\'et\'e formelle point\'ee est donn\'e par une cod\'erivation $Q:\Cal C\rightarrow\Cal C$ (c'est donc un champ de vecteurs qui s'annule en $0$). Un morphisme de vari\'et\'es formelles point\'ees est donn\'e par un morphisme de cog\`ebres. Tout morphisme $f$ de vari\'et\'es point\'ees induit un morphisme de vari\'et\'es formelles par transport des distributions de support $\{0\}$~:
$$<f_*T,\varphi>=<T,\varphi\circ f>$$
Or, par propri\'et\'e universelle des cog\`ebres cocommutatives colibres, une cod\'erivation $Q:S^+(V)\rightarrow S^+(V)$ (resp. un morphisme de cog\`ebres $\Cal F:S^+(V_1)\rightarrow S^+(V_2)$) est enti\`erement d\'etermin\'e(e) par sa composition avec la projection sur $V$ (resp. $V_2$), c'est \`a dire par une suite d'applications~:
$$Q_n:S^nV\fleche 6 V \hbox to 15mm{\hfill(resp. }\Cal F_n:S^nV_1\fleche 6 V_2\hbox{)}$$
qui sont par d\'efinition les {\sl coefficients de Taylor\/} du champ de vecteurs $Q$ ou du morphisme $\Cal F$.
\alinea{III.2. Vari\'et\'es formelles gradu\'ees point\'ees}
On fait la m\^eme construction alg\'ebrique dans la cat\'egorie des espaces vectoriels gradu\'es~: une {\sl vari\'et\'e formelle gradu\'ee point\'ee\/} est une cog\`ebre $\Cal C$ isomorphe \`a $S^+(V)$ o\`u $V$ est cette fois-ci un espace gradu\'e. Toutes les notions du \S\ III.1 s'appliquent, \`a ceci pr\`es que l'on peut consid\'erer des champs de vecteurs de diff\'erents degr\'es. Nous allons donner une formule explicite pour un champ de vecteurs ou un morphisme en fonction de ses coefficients de Taylor~:
\th{III.2.1}
Soit $i$ un entier, soient $V,V_1,V_2$ des espaces gradu\'es, et deux suites d'applications lin\'eaires $Q_n:S^nV\rightarrow V$ de degr\'e $i$, $\Cal F_n:S^nV_1\rightarrow V_2$ de degr\'e z\'ero. Alors il existe une unique cod\'erivation $Q$ de degr\'e $i$ de $S^+(V)$ et un unique morphisme $\Cal F:S^+(V_1)\rightarrow S^+(V_2)$ dont les $Q_n$ et les $\Cal F_n$ sont les coefficients de Taylor respectifs. $Q$ et $\Cal F$ sont donn\'es par les formules explicites~:
$$\eqalign{Q(x_1\ldots x_n)     &=\sum_{{\ I\amalg J=\{1,\ldots,n\}\atop I,J\not =\emptyset}}
        \varepsilon_x(I,J)\bigl(Q_{|I|}(x_I)\bigr).x_J\cr
\Cal F(x_1\ldots x_n)   &=\sum_{j\ge 1}{1\over j!}
        \sum_{{I_1\amalg\cdots\amalg I_j=\{1,\ldots,n\}\atop \uple Ij\not =\emptyset}}
\varepsilon_x(\uple Ij)\Cal F_{|I_1|}(x_{I_1})\cdots \Cal F_{|I_j|}(x_{I_j})}$$
o\`u $\varepsilon_x(\uple Ij)$ d\'esigne la signature de l'effet sur les $x_i$ impairs de la permutation-battement associ\'ee \`a la partition $(\uple Ij)$ de $\{1,\ldots,n\}$.
\dem
Supposons que tous les coefficients de Taylor de la cod\'erivation $Q$ sont nuls. En particulier $Q(x)=0$ pour tout $x\in V$. Supposons que $Q(x_1\ldots x_k)=0$ pour tout $k\le n$ Alors~:
$$\Delta Q(x_1\ldots x_{n+1})=(Q\otimes I + I\otimes Q)\Delta (x_1\ldots x_{n+1})=0,$$
compte tenu de l'hypoth\`ese de r\'ecurrence et de l'expression explicite de $\Delta(x_1\ldots x_{n+1})$~:
$$\Delta (x_1\ldots x_{n+1})=\sum_{I\amalg J=\{1,\ldots , n+1\},\ I,J\not =\emptyset}\varepsilon _x(I,J)x_I\otimes x_J$$
Donc $Q(x_1\ldots x_{n+1})\in V$, donc est nul puisque le $n+1$-\`eme coefficient de Taylor est nul. Le raisonnement est analogue dans le cas d'un morphisme, et montre qu'une cod\'erivation ou un morphisme est enti\`erement d\'etermin\'e(e) par ses coefficients de Taylor.
\ssq
Nous v\'erifions directement les formules (les v\'erifications \`a l'ordre 2 ou 3 sont laiss\'ees au lecteur \`a titre d'exercice).
\smallskip
{1. \sl Cas d'une cod\'erivation\/}~: On \'ecrit la formule explicite pour $\Delta(x_1\ldots x_n)$ en utilisant la cocommutativit\'e gradu\'ee, ce qui permet de ne retenir que la moiti\'e des partitions~:
$$\Delta(x_1\ldots x_n)=(1+\tau)\sum_{K\amalg L=\{1,\ldots ,n\},\ 1\in K, \ L\not =\emptyset}
        \varepsilon_x(K,L)x_K\otimes x_L$$
On a donc, en prenant pour $Q$ l'expression explicite du th\'eor\`eme~:
$$\eqalign{\Delta Q(x_1 & \ldots x_n)   =
                       \sum_{I\amalg J=\{1,\ldots , n\}, I,J\not=\emptyset}
        \varepsilon_x(I,J)\Delta (Q_{|I|}x_I.x_J)   \cr
                                        &=(1+\tau)
                       \sum_{I\amalg J=\{1,\ldots , n\}, I,J\not=\emptyset}
                        \ \sum_{K\amalg L=J,\ L\not =\emptyset}
\varepsilon_x(I,J)\varepsilon_{x_J}(K,L)Q_{|I|}(x_I).x_K\otimes x_L \cr
                                        &=(1+\tau)
        \ \sum_{I\amalg J\amalg K=\{1,\ldots ,n\},\ I,K\not =\emptyset}
\varepsilon_x(I,J,K)Q_{|I|}(x_I).x_J\otimes x_K     \cr}$$
Par ailleurs on a~:
$$\eqalign{(Q\otimes & I + I\otimes Q) \Delta(x_1\ldots x_n)
        =(1+\tau)\sum_{L\amalg K=\{1,\ldots , n\},\ L,K\not = \emptyset}
                \varepsilon_x(L,K)Q(x_L)\otimes x_K     \cr
        &=(1+\tau)\sum_{L\amalg K=\{1,\ldots , n\},\ L,K\not = \emptyset}
                \varepsilon_x(L,K)
        \ \sum_{I\amalg J=L,\ I,J\not =\emptyset}\varepsilon_{x_L}(I,J)
                        Q_{|I|}(x_I).x_J\otimes x_K \cr
        &=(1+\tau)
        \ \sum_{I\amalg J\amalg K=\{1,\ldots ,n\},\ ,I,J,K\not =\emptyset}
\varepsilon_x(I,J,K)Q_{|I|}(x_I).x_J\otimes x_K     \cr}$$
d'o\`u le fait que $Q$ est bien une cod\'erivation.
\smallskip
{2. \sl Cas d'un morphisme\/}~: le calcul est un peu plus compliqu\'e~: on commence par \'ecrire $\Delta$ et $\Cal F$ de mani\`ere redondante, en employant des permutations qui ne sont pas forc\'ement des battements~:
$$\eqalign{\Delta(x_1\ldots x_n)        &=\sum_{\sigma\in S_n}
                \ \sum_{r=1}^n {\varepsilon_x(\sigma)\over r!(n-r)!}
                x_{\sigma_1}\ldots x_{\sigma_r}\otimes x_{\sigma_{r+1}}\ldots x_{\sigma_n}\cr
\Cal F(x_1\ldots x_n)   &=\sum_{j\ge 1}{1\over j!}\sum_{k_1+\cdots +k_j=n}
{1\over k_1!\ldots k_j!}\sum_{\sigma\in S_n}
\varepsilon_x(\sigma)\Cal F_{k_1}(x_{\sigma_1}\ldots x_{\sigma_{k_1}})\ldots
                        \Cal F_{k_j}(x_{...}\ldots x_{\sigma_n})        \cr}$$
On v\'erifie directement l'\'egalit\'e~:
$$\Delta\Cal F(x_1\ldots x_n)=(\Cal F\otimes\Cal F)\Delta(x_1\ldots x_n)$$
L'\'ecriture par blocs~: $x_1\ldots x_n=(x_1\ldots x_{k_1})\ldots (x_{k_1+\cdots +k_{j-1}+1}\ldots x_n)$ induit par permutation des blocs un plongement du groupe de permutations $S_j$ dans $S_n$. On calcule~:
$$\eqalign{\hbox to -18pt{}\Delta\Cal F  (x_1\ldots x_n)        &=\sum_{j\ge 2}{1\over j!}
                                        \sum_{k_1+\cdots +k_j=n}
                        {1\over k_1!\ldots k_j!}\sum_{\sigma\in S_n}
\varepsilon_x(\sigma)\Delta\bigl(\Cal F_{k_1}(x_{\sigma_1}\ldots x_{\sigma_{k_1}})\ldots
                        \Cal F_{k_j}(x_{...}\ldots x_{\sigma_n})\bigr)  \cr
&=\sum_{j\ge 2}{1\over j!}
                                        \sum_{k_1+\cdots +k_j=n}
                        {1\over k_1!\ldots k_j!}\sum_{\sigma\in S_n}
\varepsilon_x(\sigma)\!\!\sum_{\tau\in S_j\subset S_n}\cr
&\hbox to 12mm{}\sum_{r=1}^{j-1}
{\varepsilon_{\sigma x}(\tau)\over r!(j-r)!}\Cal F_{k_{\tau_ 1}}(...)\ldots  \Cal F_{k_{\tau_ r}}(...)\otimes \Cal F_{k_{\tau_{r+1}}}(...)\ldots \Cal F_{k_{\tau_j}}(...)
\cr}$$
Dans le dernier membre de l'\'egalit\'e ci-dessus, chaque terme se trouve r\'ep\'et\'e autant de fois qu'il y a d'\'el\'ements dans $S_j$. On a donc~:
$$\eqalign{\Delta\Cal F  (x_1\ldots x_n)        &=
\sum_{j\ge 2}
                                        \sum_{k_1+\cdots +k_j=n}
                                {1\over k_1!\ldots k_j!}\sum_{\sigma\in S_n}
\varepsilon_x(\sigma)   \cr
&\hbox to 12mm{}\sum_{r=1}^{j-1}
{1\over r!(j-r)!}\Cal F_{k_{1}}(...)\ldots  \Cal F_{k_{r}}(...)\otimes \Cal F_{k_{{r+1}}}(...)\ldots \Cal F_{k_{j}}(...)
\cr}$$
Par ailleurs, on a~:
$$\eqalign{(\Cal F\otimes\Cal F)\Delta(x_1\ldots x_n)
&=
\sum_{\sigma\in S_n}
                \ \sum_{r=1}^{n-1} {\varepsilon_x(\sigma)\over r!(n-r)!}
                \Cal F(x_{\sigma_1}\ldots x_{\sigma_r})\otimes \Cal F(x_{\sigma_{r+1}}\ldots x_{\sigma_n})\cr
&= \sum_{\sigma\in S_n}\sum_{r=1}^{n-1}\sum_{\alpha\in S_r\times S_{n-r}\subset S_n}{\varepsilon_x(\sigma)\varepsilon_{\sigma x}(\alpha)\over r!(n-r)!}\cr
&\hbox to 12 mm{}\sum_{j,k\ge 1}{1\over j!k!}\ \sum_{{r_1+\cdots +r_j=r \atop s_1+\cdots +s_k=n-r}}{1\over r_1!\ldots r_j!s_1!\ldots s_k!}\cr
&\hbox to 12 mm{}\Cal F_{r_1}(x_{\alpha_{\sigma_1}}\ldots x_{\alpha_{\sigma_{r1}}})\ldots \Cal F_{r_j}(...)\otimes\Cal F_{s_1}(...)\ldots\Cal F_{s_k}(...x_{\alpha_{\sigma_n}})\cr}$$
Dans le dernier membre de l'\'egalit\'e ci-dessus, chaque terme se trouve r\'ep\'et\'e autant de fois qu'il y a d'\'el\'ements dans $S_r\times S_{n-r}$. Donc~:
$$\eqalign{(\Cal F\otimes\Cal F)\Delta(x_1\ldots x_n)
&=\sum_{\sigma\in S_n}
                \ \sum_{r=1}^{n-1} \varepsilon_x(\sigma)
{1\over j!k!}\ \sum_{{r_1+\cdots +r_j=r \atop s_1+\cdots +s_k=n-r}}{1\over r_1!\ldots r_j!s_1!\ldots s_k!}\cr
&\hbox to 12 mm{}\Cal F_{r_1}(x_{\sigma_1}\ldots x_{\sigma_{r_1}})\ldots \Cal F_{r_j}(...)\otimes\Cal F_{s_1}(...)\ldots\Cal F_{s_k}(...x_{\sigma_n})
\cr}$$
Posant $l=j+k$ et proc\'edant \`a la renum\'erotation $(s_1,\ldots s_k)=(r_{j+1},\ldots r_l)$ on obtient~:
$$\eqalign{(\Cal F\otimes\Cal F)\Delta(x_1\ldots x_n)
&=\sum_{\sigma\in S_n} \varepsilon_x(\sigma)\sum_{l\ge 2}
                                \sum_{r_1+\cdots +r_l=n}\sum_{k=1}^{l-1}
{1\over k!(l-k)!}{1\over r_1!\ldots r_l!}\cr
&\hbox to 12 mm{}
\Cal F_{r_1}(x_{\sigma_1}\ldots x_{\sigma_{r1}})\ldots \Cal F_{r_k}(...)\otimes\Cal F_{r_{k+1}}(...)\ldots\Cal F_{r_l}(...x_{\sigma_n})\cr
&=\Delta\Cal F(x_1\ldots x_n)\cr}$$
compte tenu du calcul pr\'ec\'edent, ce qui d\'emontre le th\'eor\`eme.
\qed
\paragraphe{IV. L$_\infty$-alg\`ebres et L$_\infty$-morphismes}
\qquad A tout espace vectoriel gradu\'e $V$ on associe (attention au d\'ecalage!) la vari\'et\'e formelle $(V[1],0)$ point\'ee, c'est \`a dire la cog\`ebre colibre sans co-unit\'e~:
$$\Cal C(V)=S^+(V[1])\tilde{\mopl{$\longmapsto$}_\Phi} \sum_{k\ge 1}(\Lambda^k V)[k]$$
o\`u $\Phi$ est l'isomorphisme d\'ecrit au \S~\ II.4.
\ssq
Un {\sl pr\'e-$L_\infty$-morphisme\/} entre deux espaces gradu\'es $V_1$ et $V_2$ est par d\'efinition un morphisme de vari\'et\'es formelles, c'est-\`a-dire un morphisme de cog\`ebres~:
$$\Cal F:\Cal C(V_1)\longrightarrow \Cal C(V_2)$$
qui est donc d\'etermin\'e par ses coefficients de Taylor $\Cal F_j$. Posant $\overline \Cal F_j=\Cal F_j\circ \Phi\inver$ on a~:
$$\displaylines{\overline\Cal F_1:V_1\longrightarrow V_2        \cr
        \overline\Cal F_2:\Lambda^2V_1\longrightarrow V_2[-1]   \cr
        \overline\Cal F_3:\Lambda^3V_1\longrightarrow V_2[-2]   \cr
                        \vdots  \cr}$$
\alinea{IV.1. Alg\`ebres de Lie homotopiques}
\qquad Par d\'efinition une {\sl $L_\infty$-alg\`ebre\/}, ou {\sl alg\`ebre de Lie homotopique\/} est une vari\'et\'e formelle gradu\'ee point\'ee du type $(\g g[1], 0)$, o\`u $\g g$ est un espace vectoriel gradu\'e, munie d'un champ de vecteurs $Q$ de degr\'e $1$ v\'erifiant l'{\sl \'equation ma\^\i tresse\/}~:
$$[Q,Q]=2Q^2=0$$
C'est-\`a-dire que $Q$ est une cod\'erivation de carr\'e nul de la cog\`ebre $\Cal C(\g g)$. Les coefficients de Taylor \hbox{$Q_k:S^k(\g g[1])\rightarrow \g g[2]$} donnent naissance aux coefficients $\overline Q_k=Q_k\circ \Phi\inver$~:
$$\displaylines{\overline Q_1:\g g\longrightarrow \g g[1]       \cr
\overline Q_2:\Lambda^2\g g\longrightarrow \g g \cr
\overline Q_3:\Lambda^3\g g\longrightarrow \g g[-1]     \cr
                        \vdots          \cr}$$
L'\'equation ma\^\i tresse se traduit par une infinit\'e de relations quadratiques entre les $\overline Q_k$, qui s'obtiennent en \'ecrivant explicitement pour tout $k$ l'\'equation~:
$$\pi Q^2(x_1\ldots x_k)=0$$
o\`u $\pi:\Cal C(\g g)\rightarrow \g g[1]$ est la projection canonique. On \'ecrit explicitement les trois premi\`eres~:
\smallskip
{\sl Premi\`ere \'equation\/}~: $Q_1^2(x)=0$ pour tout $x$ dans $\g g$. Donc $(\g g, Q_1)$ est un {\sl complexe de cocha\^\i nes\/}.
\smallskip
{\sl Deuxi\`eme \'equation\/}~: $\pi Q^2(x.y)=0$, soit~:
$$Q_2(Q_1x.y+(-1)^{|x|-1}x.Q_1y)+Q_1Q_2(x.y)=0.$$
({\sl Remarque\/}~: $|x|-1$ est bien le degr\'e de $x$ dans la cog\`ebre $\Cal C(\g g)$, \`a cause du d\'ecalage).
Traduisant cette \'egalit\'e en termes de $\overline Q_1$ et $\overline Q_2$ on obtient (cf \S\ II.4)~:
$$\alpha (\overline Q_1x,y)\overline Q_2(\overline Q_1x\wedge y)+(-1)^{|x|-1}\alpha (x,\overline Q_1y)\overline Q_2(x\wedge \overline Q_1y)+\alpha(x,y)\overline Q_1\overline Q_2(x\wedge y)=0$$
Compte tenu de l'\'egalit\'e~:
$$\alpha(x,y)=(-1)^{|x|(|y|-1)}$$
on obtient~:
$$(-1)^{|y|-1}\overline Q_2(\overline Q_1x\wedge y)-\overline Q_2(x\wedge \overline Q_1 y)+\overline Q_1\overline Q_2(x\wedge y)=0.$$
Posant $dx=(-1)^{|x|}\overline Q_1x$ et $[x,y]=\overline Q_2(x\wedge y)$ on obtient finalement~:
$$d[x,y]=[dx,y]+(-1)^{|x|}[x,dy]$$
donc $\overline Q_2$ est un crochet antisym\'etrique pour lequel $d$ est une d\'erivation.
\smallskip
{\sl Remarque\/}~: On peut garder $\overline Q_1$ comme d\'erivation sans le modifier, \`a condition d'inverser le sens du crochet, c'est-\`a-dire de poser~:
$$[x,y]=\overline Q_2(y\wedge x).$$
Nous choisirons la premi\`ere solution.
\smallskip
{\sl Troisi\`eme \'equation\/}~: $\pi Q_3(x.y.z)=0$ soit~:
$$\displaylines{Q_3\bigl(Q_1x.y.z+(-1)^{|x|-1}x.Q_1y.z+(-1)^{|x|+|y|-2}x.y.Q_1z \bigr) +Q_1Q_3(x.y.z)\hfill\cr
\hfill +Q_2\bigl(Q_2(x.y).z+(-1)^{(|y|-1)(|z|-1)}Q_2(x.z).y+(-1)^{(|x|-1)(|y|+|z|-2)}Q_2(y.z).x\bigr)=0}$$
soit~:
$$\displaylines{Q_2\bigl(Q_2(x.y).z+(-1)^{(|y|-1)(|z|-1)+(|x|-1)(|z|-1)}Q_2(z.x).y      \hfill\cr
\hfill+(-1)^{(|x|-1)(|y|+|z|)}Q_2(y.z).x\bigr)+\hbox{ \sl termes en $Q_3$ }=0\cr}$$
Or on a~:
$$\eqalign{Q_2\bigl(Q_2(x.y).z\bigr)    &=\alpha(Q_2(x.y),z)\alpha(x,y)\overline Q_2\bigl(\overline Q_2(x\wedge y)\wedge z\bigr)\cr
        &=(-1)^{(|x|+|y|)|z|}(-1)^{(|x|-1)|y|}\overline Q_2\bigl(\overline Q_2(x\wedge y)\wedge z\bigr)}$$
En reportant ceci dans l'\'equation pr\'ec\'edente et en simplifiant par $(-1)^{|x||y|+|x||z|+|y||z|}$ on obtient finalement~:
$$\permuc\toto xyz +\hbox{ \sl termes en $Q_3$ }=0$$
Autrement dit le crochet fourni par $\overline Q_2$ v\'erifie l'identit\'e de Jacobi gradu\'ee ``\`a homotopie gouvern\'ee par $\overline Q_3$ pr\`es''. En corollaire~:
\th{IV.1.1}
Une alg\`ebre de Lie diff\'erentielle gradu\'ee est la m\^eme chose qu'une $L_\infty$-alg\`ebre pour laquelle tous les coefficients de Taylor sont nuls sauf les deux premiers.
\ndem
\alinea{IV.2. L'alg\`ebre de Lie diff\'erentielle gradu\'ee des mutichamps de vecteurs}
Sur $V=T_{poly}\left(\BbbR^d\right)$, on dispose du crochet de Schouten d\'efini par~:

$$\left[\xi_1\wedge...\wedge \xi_k,\eta_1\wedge...
\wedge\eta_\ell\right]_S=\sum_{i=1}^k\sum_{j=1}^\ell (-
1)^{i+j}[\xi_i,\eta_j]\wedge\xi_1\wedge...\wedge\widehat{\xi_i}\wedge...\wedge\xi_k
\wedge\eta_1\wedge...\wedge\widehat{\eta_j}\wedge...\wedge\eta_\ell.$$

La sym\'etrisation de $Tens\left(\BbbR^d\right)$ nous permet de d\'efinir une op\'eration
$\bullet$.

Si $\alpha_1$ est un $k_1$-tenseur antisym\'etrique~:

$$\alpha_1=\alpha_1^{i_1...i_{k_1}}\partial_{i_1}\wedge\partial_{i_2}\wedge...\wedge\partial
_{i_{k_1}}\in Tens^{k_1}\left(\BbbR^d\right),$$
alors~:

$$\Phi\inver_{k_1}(\alpha_1)=\alpha_1^{i_1...i_{k_1}}\psi_{i_1}...\psi_{i_{k_1}}\in S\big({\cal
X}\left(\BbbR^d\right)[1]\big)[-k_1]$$
o\`u chaque $\psi_i=\Phi\inver_1(\partial_i)$ est une variable de degr\'e 1.

Si maintenant $\alpha_2$ est un $k_2$ tenseur antisym\'etrique, on posera~:

$$\alpha_1\bullet\alpha_2=\Phi_{k_1+k_2-1}\left(\sum_{i=1}^d{\partial
\Phi\inver_{k_1}(\alpha_1)\over \partial\psi_i}\pt {\partial\Phi\inver_{k_2}(\alpha_2)\over\partial x_i}\right)$$
en tenant compte du fait que ${\partial\over\partial\psi_i}$ est un op\'erateur de d\'erivation impair.
\lemme{IV.2.1 \rm (Calcul de $\alpha_1\bullet\alpha_2$)}
 On a~:
$$\alpha_1\bullet\alpha_2=\sum_{l=1}^{k_1}(-1)^{l-1}\alpha_1^{i_1....i_{k_1}}
\partial_l\alpha_2^{j_1...j_{k_2}}\partial_{i_1}\wedge...\wedge \widehat{\partial_{i_l}}\wedge...
\wedge\partial_{i_{k_1}}\wedge\partial_{j_1}\wedge...\wedge\partial_{j_{k_2}}$$
et

$$\left[\alpha_1,\alpha_2\right]_S=(-1)^{k_1-1}\alpha_1\bullet\alpha_2-(-1)^{k_1(k_2-
1)}\alpha_2\bullet\alpha_1.$$
\dem
On a~:
$$\eqalign{{\partial\over\partial\psi_i}\left(\alpha_1^{i_1...i_{k_1}}\psi_{i_1}\pt...\pt
\psi_{i_{k_1}}\right)&=\sum_{l=1}^{k_1}(-1)^{l-1}\alpha_1^{i_1...i_{k_1}}\psi_{i_1}\pt ...\pt
{\partial\psi_{i_l}\over \partial\psi_i}\pt...\pt\psi_{i_{k_1}}\cr
&=\sum_{l=1}^{k_1}(-1)^{l-1}\delta^i_{i_l}\alpha_1^{i_1...i_{k_1}}\psi_{i_1}\pt ...\pt
\widehat{\psi_{i_l}}\pt...\pt\psi_{i_{k_1}}.}$$
Donc~:
$$\eqalign{\sum_{i=1}^d{\partial \Phi_{k_1}\inver(\alpha_1)\over \partial\psi_i}\pt
{\partial\Phi_{k_2}\inver(\alpha_2)\over\partial x_i}&=\sum_{l=1}^{k_1}(-1)^{l-1} \alpha_1^{i_1...i_{k_1}}
\psi_{i_1}\pt ...\pt \widehat{\psi_{i_l}} \pt...\pt\psi_{i_{k_1}}\pt
\partial_{i_l}\alpha_2^{j_1...j_{k_2}} \psi_{j_1}\pt ... \pt \psi_{j_{k_2}}\cr
&=\sum_{l=1}^{k_1}(-1)^{l-1} \alpha_1^{i_1...i_{k_1}} \partial_{i_l}\alpha_2^{j_1...j_{k_2}}
\psi_{i_1}\pt ...\pt \widehat{\psi_{i_l}} \pt...\pt\psi_{i_{k_1}}\pt  \psi_{j_1}\pt ... \pt
\psi_{j_{k_2}}.}$$
D'autre part~:

$$\eqalign{\left[\alpha_1,\alpha_2\right]_S&=\big[\alpha_1^{i_1...i_{k_1}}\partial_{i_1}\wedge..
.\wedge\partial_{i_{k_1}},\alpha_2^{j_1...j_{k_2}}\partial_{j_1}\wedge...\wedge
\partial_{j_{k_2}}\big]\cr
&=\big[\alpha_1^{i_1...i_{k_1}}\partial_{i_1},\alpha_2^{j_1...j_{k_2}}\partial_{j_1}\big]\wedge
\partial_{i_2}\wedge... \wedge\partial_{i_{k_1}} \wedge\partial_{j_1}\wedge...\wedge
\partial_{j_{k_2}}+\cr
&\hskip  2cm+\sum_{l=2}^{k_2}(-
1)^{1+l}\big[\alpha_1^{i_1...i_{k_1}}\partial_{i_1},\partial_{j_l}\big]\wedge
\partial_{i_2}\wedge...\wedge\partial_{i_{k_1}}\wedge\alpha_2^{j_1...j_{k_2}}\partial_{j_1}
\wedge...\wedge\widehat{\partial_{j_l}}\wedge...\wedge\partial_{j_{k_2}}+\cr
&\hskip 2cm +\sum_{l=2}^{k_1}\big[\partial_{i_l},\alpha_2^{j_1...j_{k_2}}\partial_{j_1} \big]
\wedge\alpha_1^{i_1...i_{k_1}}\partial_{i_1}\wedge...\wedge\widehat{\partial_{i_l}} \wedge ...
\wedge\partial_{i_{k_1}}\wedge\partial_{j_2}\wedge...\wedge\partial_{j_{k_2}}\cr
&=-\sum_{l=1}^{k_2}(-1)^{l+1}\alpha_2^{j_1...j_{k_2}}\partial_{j_l}\alpha_1^{i_1...i_{k_1}}
\partial_{i_1}\wedge...\wedge\partial_{i_{k_1}}\wedge\partial_{j_1}\wedge...\wedge\widehat{
\partial_{j_l}}\wedge...\wedge\partial_{j_{k_2}}+\cr
&\hskip 2cm+\sum_{l=1}^{k_1}(-1)^{l+1}\alpha_1^{i_1...i_{k_1}}
\partial_{i_l}\alpha_2^{j_1...j_{k_2}} \partial_{j_1}\wedge\partial_{i_1}\wedge ...\wedge
\widehat{\partial_{i_l}}\wedge ...\partial_{i_{k_1}}\wedge \partial_{j_2}\wedge...
\wedge...\wedge\partial_{j_{k_2}}\cr
&=(-1)^{k_1-1}\alpha_1\bullet\alpha_2-(-1)^{(k_2-1)k_1}\alpha_2\bullet\alpha_1.}$$
\cor{IV.2.2}
L'espace gradu\'e $T_{poly}\left(\BbbR^d\right)$, muni du crochet~:

$$\left[\alpha_1,\alpha_2\right]'_S=-\left[\alpha_2,\alpha_1\right]_S$$
est aussi une alg\`ebre de Lie gradu\'ee et~:

$$\left[\alpha_1,\alpha_2\right]'_S=(-1)^{(k_1-1)k_2}\alpha_1\bullet\alpha_2+(-
1)^{k_2}\alpha_2\bullet\alpha_1.$$
\ndem
Comme $[~,~]'_S$ d\'efinit sur $T_{poly}\left(\BbbR^d\right)$ une structure d'alg\`ebre de Lie
gradu\'ee, on aura, en prenant $d=0$, une structure de $L_\infty$ alg\`ebre sur ${\cal C}
\left(T_{poly}\left(\BbbR^d\right)\right)$. Le champ de vecteurs $Q$ est caract\'eris\'e par~:

$$\eqalign{Q_1=0,\quad Q_2(\alpha_1\pt\alpha_2)&=(-1)^{(k_1-1)k_2}
\left[\alpha_1,\alpha_2\right]'_S\cr
&=\alpha_1\bullet\alpha_2+(-1)^{k_1k_2}\alpha_2\bullet\alpha_1.}$$
\alinea{IV.3. L'alg\`ebre de Lie diff\'erentielle gradu\'ee des op\'erateurs polydiff\'erentiels}
On consid\`ere l'espace vectoriel $V'=D_{poly}\left(\BbbR^d\right)$ des (combinaisons lin\'eaires d')
op\'erateurs multi\-diff\'erentiels gradu\'e par $\vert A\vert=m-1$ si $A$ est $m$-diff\'erentiel.

Sur $D_{poly}\left(\BbbR^d\right)$, l'op\'erateur de composition naturel $\circ$ s'\'ecrit~:

$$\eqalign{\left(A_1\circ A_2\right)&\left(f_1,...,f_{m_1+m_2-1}\right)=\cr
&=\sum_{j=1}^{m_1}(-1)^{(m_2-1)(j-1)} A_1\left(f_1,...,f_{j-1},A_2\left(f_j,...,f_{j+m_2-
1}\right),f_{j+m_2},...,f_{m1+m_2-1}\right).}$$
On associe \`a cette composition d'une part le crochet de Gerstenhaber~:

$$\left[A_1,A_2\right]_G=A_1\circ A_2-(-1)^{\vert A_1\vert\vert A_2\vert}A_2\circ A_1,$$
d'autre part l'op\'erateur de cobord~:

$$dA=-\left[\mu, A\right]$$
o\`u $\mu$ est la multiplication des fonctions~: $\mu(f_1,f_2)=f_1f_2$.
\smallskip
{\sl Remarque\/}~: Avec ce choix de $d$, $\big(D_{poly}\left(\BbbR^d\right),[~,~]_G,d\big)$ est une alg\`ebre de
Lie gradu\'ee diff\'erentielle, on v\'erifie en effet que $d\circ d=0$ et

$$d\left(\left[A_1,A_2\right]\right)=\big[dA_1,A_2\big]+(-1)^{\vert A_1\vert}\big[A_1,dA_2
\big].$$

L'op\'erateur de cobord de Hochschild usuel $d_H$ donn\'e par~:

$$\eqalign{\left(d_HA\right)(f_1,...,f_m)&=f_1A(f_2,...,f_m)-A(f_1f_2,f_3,...,f_m)+...+ (-
1)^mA(f_1,...,f_{m-1})f_m\cr
&=(-1)^{\vert A\vert+1}dA(f_1,...,f_m)}$$
n'est pas une d\'erivation de l'alg\`ebre de Lie gradu\'ee $\big(D_{poly}\left(\BbbR^d\right),
[~,~]_G\big)$.

Le champ de vecteurs $Q'$ sur la vari\'et\'e formelle ${\cal C}(V')$ sera donc d\'efini par~:

$$Q'_1(A)=(-1)^{\vert A\vert}dA=(-1)^{\vert A\vert+1}\left[\mu,A\right] =\left[A,\mu\right]= -
d_HA$$
et

$$\eqalign{Q'_2\left(A_1\pt A_2\right)&=(-1)^{\vert A_1\vert(\vert A_2\vert-
1)}\left[A_1,A_2\right]_G\cr
&=(-1)^{\vert A_1\vert(\vert A_2\vert-1)}A_1\circ A_2-(-1)^{\vert A_1\vert}A_2\circ A_1.}$$

\alinea{IV.4. L$_\infty$-morphismes}
Par d\'efinition un $L_\infty$-morphisme entre deux $L_\infty$-alg\`ebres $(\g g_1, Q)$ et $(\g g_2, Q')$ est un morphisme de vari\'et\'es formelles point\'ees~:
$$\Cal F:\Cal C(\g g_1)\longrightarrow \Cal C(\g g_2)$$
v\'erifiant~:
$$\Cal F Q=Q'\Cal F$$
Cette \'equation induit une infinit\'e de relations entre les coefficients de Taylor de $Q$, $Q'$ et $\Cal F$, dont nous allons examiner les deux premi\`eres~:
\smallskip
{\sl Premi\`ere \'equation\/}~: $Q_1'\Cal F_1(x)=\Cal F_1 Q_1(x)$, c'est-\`a-dire que $\Cal F_1$ est un {\sl morphisme de complexes\/}.
\smallskip
{\sl Deuxi\`eme \'equation\/}~: $\pi Q'\Cal F(x.y)=\pi \Cal F Q(x.y)$ soit~:
$$\pi Q'\bigl(\Cal F_1x.\Cal F_1 y+\Cal F_2(x.y)\bigr)=\pi \Cal F\bigl(Q_1x.y+(-1)^{(|x|-1)}x.Q_1y+Q_2(x.y)\bigr)$$
soit encore~:
\def\F{\Cal F}
$$Q'_2(\F_1x.\F_1y)+Q'_1\F_2(x.y)=\F_2\bigl(Q_1x.y+(-1)^{|x|-1}x.Q_1y\bigr)+\F_1Q_2(x.y).$$
On traduit cette derni\`ere \'egalit\'e en termes de $\overline Q_1$, $\overline Q_2$, $\overline \F_1$, etc.~:
$$\displaylines{(-1)^{|x|(|y|-1)}[\overline\F_1 x,\overline\F_1 y]
        +(-1)^{|x|+|y|-1+|x|(|y|-1)}d\overline\F_2(x\wedge y)\hfill\cr
\hfill =(-1)^{(|x|-1)(|y|-1)+|x|}\overline\F_2(dx\wedge y)
        +(-1)^{|x||y|+|x|-1+|y|}\overline\F_2(x\wedge dy)
        +(-1)^{|x|(|y|-1)}\overline\F_1([x,y])\cr}$$
soit, en multipliant par $(-1)^{|x|(|y|-1)}$~:
$$\overline\F_1([x,y])-[\overline\F_1x,\overline\F_1y]
=(-1)^{|x|+|y|-1}\bigl(d\overline\F_2(x\wedge y)-\overline\F_2(dx\wedge y)
-(-1)^{|x|}\overline\F_2(x\wedge dy)\bigr).$$
Dans le cas o\`u $\g g_1$ et $\g g_2$ sont des alg\`ebres de Lie diff\'erentielles gradu\'ees, $\F_1$ n'est donc pas forc\'ement un morphisme d'alg\`ebres de lie diff\'erentielles gradu\'ees, mais le d\'efaut est gouvern\'e par le coefficient suivant, c'est-\`a-dire $\F_2$.
\prop{IV.4.1 \rm (Equation de $L_\infty$-morphisme dans le cas des alg\`ebres de Lie diff\'erentielles gradu\'ees)}
Supposons que $\left(V,[~,~],d\right)$ et $\left(V',[~,~]',d'\right)$ soient deux alg\`ebres de
Lie gradu\'ees. Notons $\left({\cal C}(V),Q\right)$ et $\left({\cal C}(V'),Q'\right)$ les $L_\infty$
alg\`ebres correspondantes respectives. Soit ${\cal F}:{\cal C}(V)\longrightarrow{\cal C}(V')$ un
morphisme de cog\`ebre. Alors ${\cal F}$ est un $L_\infty$ morphisme si et seulement si~:
$$\eqalign{Q'_1{\cal F}_n\left(\alpha_1\pt...\pt \alpha_n\right)+&{1\over 2}\sum_{I\sqcup
J=\{1,...,n\}\atop I,J\neq\emptyset}\varepsilon_\alpha(I,J)Q'_2\left({\cal F}_{\vert
I\vert}\left(\alpha_I\right)\pt{\cal F}_{\vert J\vert}\left(\alpha_J\right)\right)=\cr
&=\sum_{k=1}^n\varepsilon_\alpha(k,1,...\hat{k},...,n){\cal F}_n\left(Q_1(\alpha_k)\pt
\alpha_1\pt...\pt\widehat{\alpha_k}\pt...\pt \alpha_n\right)+\cr
&\hskip 0.5cm+{1\over 2}\sum_{k\neq
l}\varepsilon_\alpha(k,l,1,...,\widehat{k,l},...,n){\cal F}_{n-1}\left(Q_2(\alpha_k\pt
\alpha_l)\pt\alpha_1\pt...\pt
\widehat{\alpha_k}\pt...\pt\widehat{\alpha_l}\pt...\pt\alpha_n\right)}$$
o\`u $\vert I\vert$ et $\varepsilon_\alpha(I,J)$ ont la m\^eme signification que dans le th\'eor\`eme III.2.1, et o\`u $\varepsilon_\alpha(...)$ d\'esigne le signe de Quillen de la permutation indiqu\'ee entre parenth\`eses, c'est-\`a-dire la signature de la trace sur les $\alpha_j$ impairs de cette permutation.
\ndem 
Comme pour les cod\'erivations $Q$ et les morphismes de cog\`ebres ${\cal F}$, il est facile de
voir que les applications $Q'{\cal F}$ et ${\cal F}Q$ sont uniquement d\'etermin\'ees par leur
composition avec la projection sur $V'[1]$. On d\'eduit alors l'\'equation de $L_\infty$-
morphisme
sous la forme $(Q'{\cal F})_n=({\cal F}Q)_n$ pour tout $n$. Puisqu'on est parti de deux
alg\`ebres de Lie diff\'erentielles gradu\'ees, tous les $Q_p$ et $Q'_p$ sont nuls pour $p\geq
3$.
\paragraphe{V. Quasi-isomorphismes}
\qquad Par d\'efinition un {\sl quasi-isomorphisme\/} entre deux $L_\infty$-alg\`ebres $(\g g_1,Q_1)$ et $(\g g_2,Q_2)$ est un $L_\infty$-morphisme $\Cal F$ dont le premier coefficient de Taylor $\Cal F_1:\g g_1[1]\rightarrow \g g_2[1]$ est un morphisme de complexes qui induit un isomorphisme en cohomologie (quasi-isomorphisme de complexes). Nous allons exposer la d\'emonstration du th\'eor\`eme suivant (\cite {K1} theorem 4.4)~:
\th{V.1}
Pour tout quasi-isomorphisme $\Cal F$ d'une $L_\infty$-alg\`ebre $(\g g_1,Q_1)$ vers une  $L_\infty$-alg\`ebre $(\g g_2,Q_2)$ il existe un $L_\infty$-morphisme $\Cal G$ de $(\g g_2,Q_2)$ vers $(\g g_1,Q_1)$ dont le premier coefficient de Taylor $\Cal G_1:\g g_2[1]\rightarrow \g g_1[1]$ soit un quasi-inverse pour $\Cal F_1$.
\ndem
\alinea{V.1. D\'ecomposition des $L_\infty$-alg\`ebres}
\qquad Une $L_\infty$-alg\`ebre $(\g g,Q)$ est {\sl minimale\/} si $Q_1=0$. Une $L_\infty$-alg\`ebre est {\sl lin\'eaire contractile\/} si $Q_j=0$ pour $j\ge 2$ et si la cohomologie du complexe donn\'e par $Q_1$ est triviale. On remarque que la premi\`ere notion est invariante par $L_\infty$-isomorphismes, contrairement \`a la seconde notion.
\prop{V.2}
Toute $L_\infty$-alg\`ebre $(\g g,Q)$ est $L_\infty$-isomorphe \`a la somme directe d'une $L_\infty$-alg\`ebre minimale et d'une $L_\infty$-alg\`ebre lin\'eaire contractile.
\dem
On d\'ecompose le complexe $(\g g,Q_1)$ en somme directe $(\g g',M_1)\oplus (\g g'',L_1)$ o\`u $M_1$ est une diff\'erentielle nulle et o\`u $(\g g'',L_1)$ est un complexe \`a cohomologie triviale (on n\'eglige le d\'ecalage qui n'est pas essentiel ici). Pour ce faire on note comme d'habitude $Z_k$ et $B_k$ le noyau et l'image de la diff\'erentielle en degr\'e $k$, on choisit un suppl\'ementaire $\g g'_k$ de $B_k$ dans $Z_k$, et un suppl\'ementaire $W_k$ de $Z_k$ dans $\g g_k$. Posant alors $\g g''_k=B_k\oplus W_k$ on a la d\'ecomposition cherch\'ee.
\ssq
Cette d\'ecomposition du complexe est le point de d\'epart de la d\'ecomposition de la $l_\infty$-alg\`ebre $(\g g,Q)$. La cog\`ebre associ\'ee \`a $\g g=\g g'\oplus\g g''$ s'\'ecrit~:
$$\Cal C(\g g)=\Cal C(\g g')\oplus\Cal C(\g g'')\oplus\Cal C(\g g')\otimes\Cal C(\g g'').$$
Il s'agit de construire un isomorphisme de cog\`ebres~:
$$\Cal F:\Cal C(\g g) \wt{\fleche 8} \Cal C(\g g)$$
tel que $\Cal F\circ Q=\overline Q\circ \cal F$, avec~:
$$\eqalign{&    \overline Q\restr{\Cal C(\sg g')}=M     \cr
        &\overline Q\restr{\Cal C(\sg g'')}=L           \cr
        &\overline Q\restr{\Cal C(\sg g')\otimes\Cal C(\sg g'')}=M\otimes I + I\otimes L \cr}$$
o\`u $M_1=0$, $L_j=0$ pour $j\ge 2$ et $L_1$ \`a cohomologie triviale. On pose donc pour commencer $\Cal F_1=Id:\g g\rightarrow \g g$, d'o\`u forc\'ement $\overline Q_1=Q_1=L_1$ au vu de la d\'ecomposition du complexe rappel\'ee ci-dessus. Il est tr\`es facile de voir qu'un $L_\infty$-morphisme $\Cal F$ v\'erifiant $\Cal F_1=\mop{Id}$ s'\'ecrit comme un produit infini~:
$$\Cal F=\cdots \Cal F^k\Cal F^{k-1}\cdots \Cal F^2$$
o\`u $\Cal F^k$ est le $L_\infty$-morphisme ayant l'identit\'e comme premier coefficient de Taylor, $\Cal F_k$ comme $k^{\hbox{\sevenrm i\`eme}}$ coefficient de Taylor, tous les autres coefficients \'etant nuls.
\ssq
Chercher le coefficient $\Cal F_k$ en supposant que les $\Cal F_j$ sont connus pour $j<k$, c'est donc chercher un $L_\infty$-isomorphisme $\Cal F$ ``lacunaire'' comme le $\Cal F^k$ ci-dessus, entre $(\g g'\oplus \g g'',Q)$ et $(\g g'\oplus \g g'',\overline Q)$, o\`u le champ de vecteurs impair $Q$ v\'erifie~:
$$\eqalign{&    Q_1\restr{\Cal C(\sg g')}=0     \cr
        & Q_j\bigl(\Cal C(\g g')\bigr)\subset \g g' \hbox{ pour }j\le k-1\cr
        & Q_j\bigl(\Cal C(\g g'')\bigr)=0\hbox{ pour }2\le j\le k-1\cr
        & Q_j\bigl(\Cal C(\g g')\otimes \Cal C(\g g'')\bigr)=0
                \hbox{ pour }j\le k-1   \cr}$$
et o\`u le champ de vecteurs $\overline Q$ v\'erifie les m\^emes conditions avec $k$ \`a la place de $k-1$. On supposera \'egalement que les coefficients de Taylor de $Q$ et $\overline Q$ sont les m\^emes jusqu'\`a l'ordre $k-1$ et sont nuls \`a partir de l'ordre $k+1$. Il s'agit donc simplement de trouver $\Cal F_k$ et $\overline Q_k$.
\ssq
La condition $\Cal F\circ Q=\overline Q\circ \Cal F$ s'\'ecrit, en n\'egligeant les signes provenant de la supersym\'etrie~:
$$\Cal F_k(Q_1(x_1\cdots x_k))+Q_k(x_1\cdots x_k)
        =Q_1\Cal F_k(x_1\cdots x_k)+\overline Q_k(x_1\cdots x_k)\leqno{(*)}$$
o\`u l'on a d\'esign\'e par la m\^eme lettre $Q_1$ la d\'erivation de l'alg\`ebre $S(\g g[1])$ valant $Q_1$ sur $\g g[1]$.
\smallskip
1). Si tous les $x_j,j=1\cdots k$ sont dans le noyau $Z$ de $Q_1$, l'\'equation $(*)$ se r\'eduit \`a~:
$$Q_k(x_1\cdots x_k)=Q_1\Cal F_k(x_1\cdots x_k)+\overline Q_k(x_1\cdots x_k).
        \leqno{(*)_1}$$
On choisit donc $\overline Q_k(x_1\cdots x_k)$ comme \'etant la projection de $Q_k(x_1\cdots x_k)$ sur le suppl\'emen\-taire $W\oplus \g g'$ de $B$ dans $\g g$. Ceci permet de d\'efinir $\Cal F_k(x_1\cdots x_k)$ \`a un \'el\'ement $z$ de $Z$ pr\`es.
\ssq
On utilise alors l'\'equation ma\^\i tresse $[Q,Q]=0$, qui permet de montrer, par r\'ecurrence sur $k$, que $Q_k(x_1\cdots x_k)$ appartient \`a $Z$. On en d\'eduit que $\overline Q_k(x_1\cdots x_k)$ appartient bien \`a $\g g'$.
\ssq
De plus si $x_1=Q_1y_1\in B$, l'\'equation $[Q,Q]=0$ s'\'ecrit (toujours en n\'egligeant les probl\`emes de signes)~:
$$Q_k(Q_1y_1.x_2\cdots x_k)+\hbox{ termes interm\'ediaires }+Q_1Q_k(y_1.x_2\cdots x_k)=0.$$
Les termes interm\'ediaires sont une somme de termes du type~:
$$Q_j(\cdots Q_l(\cdots)\cdots), j,l<k.$$
L'\'el\'ement $Q_1y_1$ se trouve dans une parenth\`ese int\'erieure ou dans la parenth\`ese ext\'erieure. Dans les deux cas l'hypoth\`ese de d\'epart sur $Q$ entra\^\i ne l'annulation de ce terme. On a donc~:
$$Q_k(Q_1y_1.x_2\cdots x_k)=-Q_1Q_k(y_1.x_2\cdots x_k),$$
ce qui montre que $\overline Q_k(Q_1y_1.x_2\cdots x_k)=0$.
\medskip
2). Soit $s\in S^k(\g g[1])$, avec $k\ge 2$. On dit que $x$ est de type $j, 0\le j\le k$, si $x$ s'\'ecrit $x_1\cdots x_k$ avec $\uple xj\in W$ et $x_{j+1},\ldots ,x_k\in Z$. Nous allons d\'eterminer $\Cal F_k(x)$ par r\'ecurrence (finie) sur le type de $x$, le type $0$ ayant \'et\'e trait\'e au 1). On remarque que $\overline Q_k(x)=0$ si le type de $x$ est non nul.
\ssq
L'\'equation $[Q,Q]=0$ s'\'ecrit~:
$$Q_k(Q_1(x_1\cdots x_k))+Q_1Q_k(x_1\cdots x_k)=0,$$
les termes interm\'ediaires s'annulant pour la m\^eme raison que dans le 1). On a donc~:
$$Q_1Q_k(x)+Q_kQ_1(x)=0 \leqno{(M)}$$
pour tout $x\in S^k(\g g[1])$.
\ssq
Soit $r\ge 1$. Supposons que $\Cal F_k(x)$ soit d\'etermin\'e pour tout $x$ de type $j\le r-2$, et d\'etermin\'e \`a un $z\in Z$ pr\`es pour tout $x$ de type $r-1$. Soit alors $x$ de type $r$. On veut d\'eterminer $\Cal F_k(x)$ \`a un \'el\'ement $z'\in Z$ pr\`es et pr\'eciser $\Cal F_k(y)$ pour tous les $y$ de type $r-1$.
\ssq
L'\'equation $(*)$ appliqu\'ee \`a $Q_1(x)$ s'\'ecrit~:
$$\Cal F_kQ_1^2(x)+Q_kQ_1(x)=Q_1\Cal F_k Q_1(x)$$
le terme $\overline Q_kQ_1(x)$ \'etant nul. En reportant (M) dans cette \'equation on a donc~:
$$Q_1\Cal F_k Q_1(x)+Q_1Q_k(x)=0,$$
d'o\`u~:
$$\Cal F_kQ_1(x)+Q_k(x)\in Z.$$
Comme $\Cal F_kQ_1(x)$ est d\'etermin\'e \`a un \'el\'ement arbitraire de $Z$ pr\`es, on peut s'arranger pour que~:
$$\Cal F_k Q_1(x)+Q_k(x)=b(x)$$
o\`u $b(x)$ appartient \`a $B$. L'\'equation $(*)$ appliqu\'ee \`a $x$ s'\'ecrivant~:
$$\Cal F_kQ_1(x)+Q_k(x)=Q_1\Cal F_k(x)$$
le choix d'un $b(x)$ nous permet de choisir $\Cal F_k(x)$ \`a un \'el\'ement $z'\in Z$ pr\`es. Le $b(x)$ doit ob\'eir \`a la contrainte suivante~: si $Q_1(x)=0$, alors $b(x)=Q_k(x)$. Supposons que $x=Q_1y$ o\`u $y$ est de type $r+1$. Alors, compte tenu de (M) la contrainte sur $b$ s'\'ecrit~:
$$b(x)=-Q_1Q_k(y).$$
Ayant choisi un $b(x)$ pour tout $x$ de type $r$ satisfaisant \`a la contrainte ci-dessus, on peut alors choisir $\Cal F_k(x)$ \`a un \'el\'ement $z'\in Z$ pr\`es. Il reste donc simplement \`a d\'emontrer le lemme ci-dessous~:
\lemme{V.3}
Soit $x$ de type $r\ge 1$. Alors si $Q_1x=0$ il existe un $y$ de type $r+1$ tel que $x=Q_1y$.
\dem
On consid\`ere l'application $\delta :\g g\rightarrow \g g$ de degr\'e $-1$ d\'efinie par $\delta(x)=0$ pour $x\in \g g'\oplus W$, et $\delta(Q_1x)=x$ pour tout $x$ dans $\g g$. On a alors~:
$$Q_1\delta+\delta Q_1=\mop{Id}-p,$$
o\`u $p$ est la projection sur $\g g'$ parall\`element \`a $\g g''$ (autrement dit $\delta$ est une homotopie entre les deux endomorphismes de complexes $\mop{Id}$ et $p$).
\ssq
Le lemme V.3 est un corollaire du r\'esultat suivant, d\^u \`a Quillen [Q appendix B]~:
\prop{V.4}
1). La d\'erivation $Q_1$ de l'alg\`ebre sym\'etrique $S(\g g)$ v\'erifie~:
$$Q_1^2=0.$$
2). La cohomologie du complexe $(S(\g g),Q_1)$ est isomorphe \`a $S(\g g')$, et un suppl\'ementaire de l'image de $Q_1$ dans le noyau de $Q_1$ est donn\'e par $S(\g g')\otimes 1$ moyennant l'identification~: $S(\g g)=S(\g g')\otimes S(\g g'')$.
\dem

1). Comme $Q_1$ est impaire, $Q_1^2=\frac 12 [Q_1,Q_1]$ est encore une d\'erivation de $S(\g g)$. Comme $Q_1^2\restr{\sg g}=0$ cette d\'erivation est nulle.
\smallskip
2). On a~: $Q_1(v'v'')=v'Q_1(v'')$ pour $v'\in S(\g g')$ et $v''\in S(\g g'')$. On est ramen\'e au cas o\`u la cohomologie de $\g g$ est triviale. On prolonge alors l'homotopie $\delta$ ci-dessus en une d\'erivation de $S(\g g)$. On pose alors~:
$$E=[Q_1,\delta]=Q_1\delta+\delta Q_1.$$
$E$ est une d\'erivation telle que $E\restr {\sg g}=\mop{Id}$. On en d\'eduit~:
$$E(x)=kx$$
pour tout $x\in S^k(\g g)$. Si maintenant $x$ appartient \`a $S^k(\g g)$ et $Q_1 x=0$, alors $Ex=Q_1\delta x=kx$. Si $k\ge 1$ on a donc~:
$$x=Q_1(\frac 1k \delta x).$$
La cohomologie de $S(\g g)$ est donc r\'eduite au corps de base, qui est $S(\{0\})$.
\qed
{\sl Fin de la d\'emonstration du lemme V.3\/}~: le complexe $V=S(\g g)$ admet \`a son tour une d\'ecomposition~:
$$V=V'\oplus V''$$
avec $V'=S(\g g')\otimes 1$. L'image de $Q_1$ dans $S(\g g)$ est l'id\'eal engendr\'e par $B=Q_1(\g g)$. On peut donc choisir pour $V''$ l'id\'eal engendr\'e par $\g g''$. Le lemme provient alors du fait que tout \'el\'ement de type $r\ge 1$ appartient \`a cet id\'eal, sur lequel la cohomologie est triviale.
\qed
\alinea{V.2. D\'emonstration du th\'eor\`eme V.1}
On se donne deux $L_\infty$-alg\`ebres $(\g g_1,Q_1)$ et $(\g g_2,Q_2)$ et un quasi-isomorphisme $\Cal F$ de $(\g g_1,Q_1)$ vers $(\g g_2,Q_2)$. Appliquant la proposition $V.2$ \`a ces deux $L_\infty$-alg\`ebres on a le diagramme suivant, dans lequel toutes les fl\`eches sont des quasi-isomorphismes~:
$$\Cal C({\g g'}_1)\inj 7_i \Cal C({\g g'}_1\oplus {\g g''}_1)
                \widetilde{\fleche 7} \Cal C(\g g_1)
                \fleche 7_{\Cal F} \Cal C(\g g_2)
                \widetilde{\fleche 7} \Cal C({\g g'}_2\oplus {\g g''}_2)
                \surj 7_p \Cal C({\g g'}_2).$$
On a ainsi construit un quasi-isomorphisme $\Cal F'$ entre deux $L_\infty$-alg\`ebres minimales. Son premier coefficient $\Cal F'_1:\g g'_1\rightarrow \g g'_2$ \'etant inversible, $\Cal F'$ lui-m\^eme est inversible. L'ajout du quasi-isomorphisme ${\Cal F'}\inver$ dans le diagramme ci-dessus permet alors la construction d'un quasi-isomorphisme~:
$$\Cal G:\Cal C(\g g_2)\longrightarrow \Cal C(\g g_1)$$
qui est un quasi-inverse pour $\Cal F$.
\qed
\paragraphe{VI. La formalit\'e de Kontsevich}
Un $L_\infty$ morphisme entre $T_{poly}\left(\BbbR^d\right)$ et
$D_{poly}\left(\BbbR^d\right)$ qui soit aussi un quasi-isomorphisme c'est \`a dire un
isomorphisme en cohomologie est une {\sl formalit\'e}.

M. Kontsevich a propos\'e dans \cite {K1} une formalit\'e ${\cal U}$ explicite. Pr\'ecis\'ement, les
applications ${\cal U}_n$ sont donn\'es par~:

$${\cal U}_n=\sum_{m\geq 0}\sum_{\vec{\Gamma}\in G_{n,m}}w_{\vec{\Gamma}}{\cal
B}_{\vec{\Gamma}}$$
o\`u $G_{n,m}$ est l'ensemble des graphes orient\'es admissibles \`a $n$ sommets a\'eriens
$p_1$,...,$p_n$ et $m$ sommets terrestres $q_1$,...,$q_m$~: de chaque sommet a\'erien est issu
$k_1$,..., $k_n$ fl\`eches aboutissant soit \`a un autre sommet a\'erien soit \`a un sommet
terrestre. On ordonne les sommets a\'eriens et terrestres du graphe et on oriente le graphe en
ordonnant les fl\`eches de fa\c con compatible avec cet ordre, les fl\`eches issues du sommet
$p_j$ ont les num\'eros $k_1+...+k_{j-1}+1$,..., $k_1+...+k_j$. On les note~:

$$Star(p_j)=\{\overrightarrow{p_ja_1},...,\overrightarrow{p_ja_{k_j}}\}\qquad
\overrightarrow{v}_{k_1+...+k_{j-1}+i}=\overrightarrow{p_ja_i}.$$
Si $\vec{\Gamma}$ est un graphe orient\'e, son {\sl poids\/} $w_{\vec{\Gamma}}$ est
par d\'efinition l'int\'egrale sur l'espace de configuration $C^+_{\{p_1,...,p_n\},\{q_1,...,q_m\}}$
de la forme~:

$$\omega_{\vec{\Gamma}}={1\over (2\pi)^{\sum k_i}k_1!...k_n!}d\Phi_{\overrightarrow{v}_1}
\wedge ...\wedge d\Phi_{\overrightarrow{v}_{k_1+...+k_n}}\quad\hbox{o\`u}\quad 
\Phi_{\overrightarrow{p_ja}}=Arg\left({a-p_j\over a-\overline{p_j}}\right).$$

Enfin ${\cal B}_{\vec{\Gamma}}$ est un op\'erateur $m$-diff\'erentiel, nul sur
$\alpha_1\pt...\pt\alpha_n$ sauf si $\alpha_1$ est un $k_1$-tenseur, $\alpha_2$ un $k_2$-
tenseur,...,  $\alpha_n$ un $k_n$-tenseur, auquel cas, on a~:

$${\cal B}_{\vec{\Gamma}}(\alpha_1\pt...\pt\alpha_n)(f_1,f_2,...,f_m)=\sum D_{p_1}
\alpha_1^{i_1i_2...i_{k_1}}...D_{p_n}\alpha_n^{i_{k_1+...k_{n-
1}+1}...i_{k_1+...+k_n}}D_{q_1}f_1...D_{q_m}f_m$$
si $D_a$ est l'op\'erateur~:

$$D_a=\prod_{l,\overrightarrow{v_l}=\overrightarrow{.a}}\partial_{i_l}$$
et si la somme est \'etendue \`a tous les indices $i_j$ r\'ep\'et\'es. On notera aussi 

$${\cal U}_n=\sum{\cal U}_{(k_1,k_2,...,k_n)}=\sum{\cal U}_{k_{\{1,...,n\}}}.$$
Maintenant, si on change l'ordre des fl\`eches issues d'un sommet $p_j$, le produit
$w_{\vec{\Gamma}}{\cal
B}_{\vec{\Gamma}}$ ne change pas. On prend la convention suivante~: si $\vec{\Gamma}$ est
un graphe orient\'e de fa\c con non
compatible, on pose~:

$$B_{\vec{\Gamma}}=\varepsilon(\sigma)B_{\vec{\Gamma}^\sigma}$$
o\`u $\sigma$ est n'importe quelle permutation des fl\`eches de $\vec{\Gamma}$ qui le
transforme en un graphe $\vec{\Gamma}^\sigma$
orient\'e de fa\c con compatible. Avec cette convention, on aura~:

$${\cal U}_n=\sum_{m\geq 0}\sum_{\vec{\Gamma}\in G'_{n,m}}w'_{\vec{\Gamma}}{\cal
B}_{\vec{\Gamma}'}$$
o\`u $G'_{n,m}$ est l'ensemble de tous les graphes orient\'es de fa\c con compatible ou non et
$w'_{\vec{\Gamma}}$ est
l'int\'egrale de la forme~:

$$\omega'_{\vec{\Gamma}}={1\over (2\pi)^{\sum k_i}(\sum k_i)!}d\Phi_{\overrightarrow{v}_1}
\wedge ...\wedge d\Phi_{\overrightarrow{v}_{k_1+...+k_n}}\quad\hbox{o\`u}\quad 
\Phi_{\overrightarrow{p_ja}}=Arg\left({a-p_j\over a-\overline{p_j}}\right).$$

Nous allons v\'erifier dans la suite que nos choix de signes sont coh\'erents.
\th{VI.1 \rm(M.Kontsevich)}
L'application formelle ${\cal U}$ est une formalit\'e. En particulier c'est un $L_\infty$-morphisme.
\dem
Puisque $Q_1=0$, l'\'equation de formalit\'e s'\'ecrit~:

$$\eqalign{0=&Q'_1\left({\cal U}_{k_{\{1,...,n\}}}(\alpha_1\pt...\pt\alpha_n)\right)+\cr
&+{1\over 2}\sum_{I\sqcup J=\{1,...,n\}\atop I,J\neq\emptyset}\varepsilon_\alpha(I,J)
Q'_2\left({\cal U}_{k_I}(\alpha_I)\circ{\cal U}_{k_J}(\alpha_J)\right)-\cr
&-{1\over 2}\sum_{i\neq j}\varepsilon_\alpha (i,j,1,...,\widehat{i,j},...,n){\cal
U}_{((k_i+k_j-1), k_1,..., \widehat{k_i} ,... ,\widehat{k_j},...,k_n)}\left(Q_2(\alpha_i\pt \alpha_j)
\pt\alpha_1 \pt...\pt \widehat{\alpha_i} \pt...\pt \widehat{\alpha_j}\pt...\pt\alpha_n\right).}$$
Remarquons maintenant que pour que $w_{\vec{\Gamma}}$ ne soit pas nul, il faut que le degr\'e de la forme
$\omega_{\vec{\Gamma}}$ soit \'egal \`a la dimension de l'espace de configuration $C^+_{\{p_1,...,p_n\};\{q_1,...,q_m\}}$ sur lequel
on int\`egre, c'est \`a dire~:

$$\sum k_i=2n+m-2.$$
Dans ce cas,

$$(-1)^m=(-1)^{\left\vert{\cal U}_{k_{\{1,\ldots,n\}}}(\cdots)\right\vert+1}=(-1)^{\sum k_i}=(-1)^{\vert
k_{\{\}1,\ldots n}\vert}.$$
Donc notre \'equation devient~:

$$\eqalign{(1)&0={\cal U}_{k_{\{1,...,n\}}}(\alpha_{\{1,...,n\}})\circ\mu-(-1)^{\sum k_i-1}\mu
\circ{\cal U}_{k_{\{1,...,n\}}}(\alpha_{\{1,...,n\}})+\cr
(2)&+{1\over 2}\sum_{I\sqcup J=\{1,...,n\}\atop I,J\neq\emptyset}\varepsilon_\alpha(I,J) (-
1)^{(\vert k_I\vert-1)\vert k_J\vert}{\cal U}_{k_I}(\alpha_I)\circ{\cal U}_{k_J}(\alpha_J)+\cr
(3)&+{1\over 2}\sum_{I\sqcup J=\{1,...,n\}\atop I,J\neq\emptyset}\varepsilon_\alpha(I,J) (-
1)^{\vert k_I\vert}{\cal U}_{k_J}(\alpha_J)\circ{\cal U}_{k_I}(\alpha_I)-\cr
(4)&-{1\over 2}\sum_{i\neq j}\varepsilon_\alpha (i,j,1,...,\widehat{i,j},...,n){\cal
U}_{((k_i+k_j-1), k_1,..., \widehat{k_i} ,... ,\widehat{k_j},...,k_n)}\left((\alpha_i\bullet \alpha_j)
\pt\alpha_1 \pt...\pt \widehat{\alpha_i} \pt...\pt \widehat{\alpha_j}\pt...\pt\alpha_n\right)-\cr
(5)&-{1\over 2}\sum_{i\neq j}\varepsilon_\alpha (i,j,1,...,\widehat{i,j},...,n){\cal
U}_{((k_i+k_j-1),..., \widehat{k_i} ,... ,\widehat{k_j},...,k_n)}\left((-1)^{k_ik_j}(\alpha_j\bullet
\alpha_i) \pt...\pt \widehat{\alpha_i} \pt...\pt \widehat{\alpha_j}\pt...\pt\alpha_n\right).}$$

Montrons que $(2)=(3)$. En fait~:

$$\varepsilon_\alpha(I,J)=\varepsilon_\alpha(J,I)(-1)^{\vert k_I\vert\vert k_J\vert}$$
car le nombre de $i$ de $I$ tel que $k_i-2$ soit impair est congru modulo \`a 2 \`a $\vert
k_I\vert=\sum k_i$. Donc~:

$$(3)={1\over 2}\sum_{I\sqcup J=\{1,...,n\}\atop I,J\neq\emptyset}\varepsilon_\alpha(J,I) (-
1)^{\vert k_I\vert\vert k_J\vert+\vert k_I\vert}{\cal U}_{k_J}(\alpha_J)\circ{\cal
U}_{k_I}(\alpha_I)$$
en changeant les r\^oles de $I$ et $J$~:

$$(3)={1\over 2}\sum_{I\sqcup J=\{1,...,n\}\atop I,J\neq\emptyset}\varepsilon_\alpha(I,J) (-
1)^{\vert k_J\vert(\vert k_I\vert-1)}{\cal U}_{k_I}(\alpha_I)\circ{\cal U}_{k_J}(\alpha_J)=(2).$$
De m\^eme $(5)=(4)$~:

$$\eqalign{&(4)=-{1\over 2}\sum_{i\neq j}\varepsilon_\alpha (i,j,1,...,\widehat{i,j},...,n){\cal
U}_{((k_i+k_j-1), k_1,..., \widehat{k_i} ,... ,\widehat{k_j},...,k_n)}\left((\alpha_i\bullet \alpha_j)
\pt\alpha_1 \pt...\pt \widehat{\alpha_i} \pt...\pt \widehat{\alpha_j}\pt...\pt\alpha_n\right)\cr
&=-{1\over 2}\sum_{i\neq j}\varepsilon_\alpha (j,i,1,...,\widehat{i,j},...,n)(-1)^{k_ik_j}{\cal
U}_{((k_i+k_j-1), k_1,..., \widehat{k_i} ,... ,\widehat{k_j},...,k_n)}\left((\alpha_i\bullet \alpha_j)
\pt\alpha_1 \pt...\pt \widehat{\alpha_i} \pt...\pt \widehat{\alpha_j}\pt...\pt\alpha_n\right)\cr
&=-{1\over 2}\sum_{i\neq j}\varepsilon_\alpha (i,j,1,...,\widehat{i,j},...,n){\cal
U}_{((k_i+k_j-1), k_1,..., \widehat{k_j} ,... ,\widehat{k_i},...,k_n)}\left((-1)^{k_ik_j}(\alpha_j\bullet
\alpha_i)
\pt\alpha_1 \pt...\pt \widehat{\alpha_j} \pt...\pt \widehat{\alpha_i}\pt...\pt\alpha_n\right)\cr
&=(5).}$$
Posons enfin $\mu={\cal U}_\emptyset$. Alors $(1)$ s'\'ecrit~:

$$(1)=(-1)^{(\vert k_{\{1,...,n\}}\vert-1).0}\ {\cal
U}_{k_{\{1,...,n\}}}\left(\alpha_{\{1,...,n\}}\right)\circ{\cal U}_\emptyset+
(-1)^{(0-1)\vert k_{\{1,...,n\}}\vert}{\cal U}_\emptyset\circ{\cal
U}_{k_{\{1,...,n\}}}\left(\alpha_{k_{\{1,...,n\}}}\right).$$
Comme $\varepsilon_\alpha(\{1,...,n\},\emptyset)=\varepsilon_\alpha(\emptyset,\{1,...,n\})=1$,
l'\'equation de formalit\'e devient~:

$$\eqalign{\sum_{I\sqcup J=\{1,...,n\}}&\varepsilon_\alpha(I,J) (-1)^{(\vert k_I\vert-1)\vert
k_J\vert}{\cal
U}_{k_I}(\alpha_I)\circ{\cal U}_{k_J}(\alpha_J)-\sum_{i\neq j}\varepsilon_\alpha (i,j,1,...,\widehat{i,j},...,n)\cr
&\hskip 2cm{\cal
U}_{((k_i+k_j-1), k_1,..., \widehat{k_i} ,... ,\widehat{k_j},...,k_n)}\left((\alpha_i\bullet \alpha_j)
\pt\alpha_1 \pt...\pt \widehat{\alpha_i} \pt...\pt
\widehat{\alpha_j}\pt...\pt\alpha_n\right)=0.}$$
Si on remplace les ${\cal U}_I(\alpha_I)$ par les
$\sum_{\vec{\Gamma}}w'_{\vec{\Gamma}}{\cal B}'_{\vec{\Gamma}}(\alpha_I)$ et qu'on
d\'eveloppe tout, on obtient une somme d'op\'erateurs multi-diff\'erentiels de la forme~:

$$\sum_{\vec{\Gamma'}}c_{\vec{\Gamma'}}{\cal
B}'_{\vec{\Gamma'}}(\alpha_1\pt...\pt\alpha_n)$$
o\`u $\vec{\Gamma'}$ est un graphe \`a $n$ sommets a\'eriens, $m$ sommets terrestres ayant
$2n+m-3$ fl\`eches. Si on se donne
$\vec{\Gamma'}$ orient\'e et une face $F$ de codimension 1 de $\partial
C^+_{\{p_1,...,p_n\};\{q_1,...,q_m\}}$, on associe \`a ce
couple $(\Gamma', F)$ au plus un terme de l'\'equation de formalit\'e. Plus pr\'ecis\'ement~:

\vskip 0.5cm
{\bf Cas 1}~: si

$$\eqalign{F&=\partial_{\{p_{i_1},...,p_{i_{n_1}}\};\{q_{l+1},...,q_{l+m_1}\}}
C^+_{\{p_1,...,p_n\};\{q_1,...,q_m\}}\cr
&=C^+_{\{p_{i_1},...,p_{i_{n_1}}\};\{q_{l+1},...,q_{l+m_1}\}}\times
C^+_{\{p_1,...,p_n\}\setminus\{p_{i_1},...,p_{i_{n_1}}\};\{q_1,...,q_l,q,q_{l+m_1+1},...,q_m\}},}$$
que l'on notera~:
$$\partial_{S,S'}C^+_{A,B}=C^+_{S,S'}\times C^+_{A\setminus S,\ B\setminus S'\sqcup \{q\}},$$
on associe au couple $(\vec{\Gamma}',F)$ l'unique terme :

$${\cal B}'_{\vec{\Gamma}',F}(\alpha_1,..,\alpha_n)(f_1,...,f_m)={\cal
B}'_{\vec{\Gamma}_2}\left(\alpha_{j_1}\pt...\pt\alpha_{j_{n_2}}\right) \left(f_1,...,f_l,{\cal
B}'_{\vec{\Gamma}_1}\left(\alpha_{i_1}\pt...\pt\alpha_{i_{n_1}}\right)
\left(f_{l+1},...,f_{l+m_1}\right),f_{l+m_1+1},...f_m\right)$$
o\`u $\vec{\Gamma}_1$ est la restriction \`a $\{p_{i_1},...,p_{i_{n_1}}\}\cup
\{q_{l+1},...,q_{l+m_1}\}$ (avec son ordre),
$\vec{\Gamma}_2$ est le graphe obtenu en collapsant les points $p_{i_1}$,..., $p_{i_{n_1}}$ et
$q_{l+1}$,..., $q_{l+m_1}$ en $q$, on
a pos\'e $\{1,...,n\}\setminus\{i_1,...,i_{n_1}\}=\{j_1<j_2<...<j_{n_2}\}$. On note
$c_{\vec{\Gamma}',F}$ le coefficient de cet op\'erateur.

Remarquons que l'application
$(\vec{\Gamma}',F)\mapsto(\vec{\Gamma}_1,\vec{\Gamma}_2)$ est dans ce cas surjective
mais pas injective. Si on se donne le couple $(\vec{\Gamma}_1,\vec{\Gamma}_2)$, la face
$F$ est bien d\'etermin\'ee mais $\vec{\Gamma'}$ n'est pas unique~: il y a d'abord la r\'epartition
des fl\`eches allant d'un sommet de $\vec{\Gamma}_2$ vers un sommet de
$\vec{\Gamma}_1$ (application de la r\`egle de Leibniz) chaque r\'epartition correspond \`a un
graphe $\Gamma_2$ diff\'erent. Si cette r\'epartition est donn\'ee, il faut encore fixer l'ordre des
fl\`eches de $\Gamma'$. Le nombre de choix est bien s\^ur le quotient du nombre d'orientations
possibles pour $\Gamma'$ par celui des orientations possibles de $\Gamma_1$ et
$\Gamma_2$~:

$$\hbox{Nombre d'orientations de $\Gamma'$}=
{\left(k_{i_1}+...+k_{i_{n_1}}\right)!\left(k_{j_1}+...+k_{j_{n_2}}\right)!\over \left(\sum
k_i\right)!}= {\vert k_I\vert!\vert k_J\vert!\over \left\vert k_{\{1,...,n\}}\right\vert !}.$$

\vskip 0.5cm
{\bf Cas 2}~: si

$$F=\partial_{\{p_i,p_j\}} C^+_{\{p_1,...,p_n\};\{q_1,...,q_m\}}=C_{\{p_i,p_j\}}\times
C^+_{\{p,p_1,...,\widehat{p_i},...,\widehat{p_j},...,p_n\};\{q_1,...,q_m\}},$$
que l'on notera~:
$$\partial_S C^+_{A,B}=C_S\times C^+_{A\setminus S\sqcup \{p\},\ B},$$
A $(\vec{\Gamma}',F)$, si la fl\`eche $\overrightarrow{p_ip_j}$ est une des fl\`eches de
$\Gamma'$, on associe l'unique terme :

$${\cal B}'_{\vec{\Gamma}',F}(\alpha_1,..,\alpha_n)(f_1,...,f_m)={\cal B}_{\vec{\Gamma}_2}
\left((\alpha_i\bullet\alpha_j)\pt\alpha_1\pt...\pt\widehat{\alpha_i}\pt...\pt\widehat{\alpha_j}
\pt...\pt\alpha_n\right)$$
o\`u $\vec{\Gamma}_2$ est le graphe obtenu en collapsant les sommets $p_i$ et $p_j$ du
graphe $\Gamma'$ sur le point $p$ et en \'eliminant la fl\`eche $\overrightarrow{p_ip_j}$. Si
cette fl\`eche n'existe pas dans $\Gamma'$, on associe l'op\'erateur nul \`a
$(\vec{\Gamma}',F)$. On note $c_{\vec{\Gamma}',F}$ le coefficient de cet op\'erateur.

Dans ce cas, on consid\'erera l'application $(\vec{\Gamma}',F)\mapsto
(\vec{\Gamma}_1,\vec{\Gamma}_2)$ o\'u $\vec{\Gamma}_1$ est le graphe trac\'e dans
$C_{\{p_i,p_j\}}$ \`a une seule fl\`eche~: la fl\`eche $\overrightarrow{p_ip_j}$. A part le cas 0,
l'image r\'eciproque d'un couple $(\vec{\Gamma}_1,\vec{\Gamma}_2)$ contient exactement~:

$$\hbox{Nombre d'orientations de $\Gamma'$}= {\left((k_i+k_j-1)+k_1+...+\widehat{k_i}+...+
\widehat{k_j}+...+k_n\right)!\over \left(\sum k_i\right)!}= {\left(\vert k_{\{1,...,n\}}\vert-1\right)!
\over \left\vert k_{\{1,...,n\}}\right\vert !}.$$

\vskip 0.5cm
\noindent
{\bf Cas 3}~: si

$$F=\partial_S C^+_{A,B}=C_S\times C^+_{A\setminus S\sqcup \{p\},\ B}$$
avec $|S|\ge 3$, dans ce cas aucun terme de l'\'equation de formalit\'e n'est associ\'e \`a $(\vec{\Gamma}',F)$.
On pose donc $c_{\vec{\Gamma}',F}=0$.

Pour chaque $\vec{\Gamma}'$, on d\'efinit sur $C^+_{\{p_1,...,p_n\};\{q_1,...,q_m\}}$ la
forme~:

$$\omega'_{\vec{\Gamma}'}={1\over (2\pi)^{\left\vert k_{\{1,...,n\}}\right\vert}\left\vert
k_{\{1,...,n\}}\right\vert !}d\Phi_{\overrightarrow{v}_1}
\wedge ...\wedge d\Phi_{\overrightarrow{v}_{k_1+...+k_n}}.$$
Montrons qu'avec toutes ces notations, l'\'equation de formalit\'e s'\'ecrit :

$$\eqalign{0&=\sum_{\vec{\Gamma}'\in G'_{n,m}} \left[\sum_{F\in
\partial C^+_{A,B}}\int_{\vec{F}} \omega'_{\vec{\Gamma}'}\right]
{\cal B}'_{\vec{\Gamma}'}\left(\alpha_1\pt...\pt\alpha_n\right)\cr
&=\sum_{\vec{\Gamma}'\in G'_{n,m}}
\left[\int_{C^+_{A,B}}d\omega'_{\vec {\Gamma}}\right]{\cal
B}'_{\vec{\Gamma}'}\left(\alpha_1\pt...\pt\alpha_n\right).}$$
Le r\'esultat est donc une simple cons\'equence du th\'eor\`eme de Stokes sur la vari\'et\'e \`a
coins $C^+_{A,B}$ et pour les formes ferm\'ees
$\omega'_{\vec{\Gamma}'}$.

Comparons donc terme par terme chaque coefficient $c_{\vec{\Gamma}',F}$ et l'int\'egrale sur
la face orient\'ee $\vec{F}$ de la forme $\omega'_{\vec{\Gamma}'}$.

\vskip 0.5cm
{\bf Cas 1}~: Avec nos notations, le coefficient $c_{\vec{\Gamma}',F}$ est~:

$$\eqalign{c_{\vec{\Gamma}',F}&=\varepsilon_\alpha(J,I)(-1)^{(\vert k_J\vert-1) \vert
k_I\vert}{\vert k_I\vert !\vert k_J\vert !\over \left\vert k_{\{1,...,n\}}\right\vert !} (-1)^{l(m_1-1)}
\cr
&\hskip 2cm \int_{C^+_{S,S'}}
\omega'_{\vec{\Gamma}_1} \int_{C^+_{A\setminus S,\ B\setminus S'\sqcup \{q\}}} \omega'_{\vec{\Gamma}_2}.}$$
(Le signe $(-1)^{l(m_1-1)}$ provient du d\'eveloppement de l'op\'eration $\circ$). On rappelle que $S'$ est le segment ${q_{l+1},\ldots ,q_{l+m_1}}$ et que $q$ remplace $S'$ dans $B\setminus S'\sqcup \{q\}$. D'autre part, on
a pour la forme

$$\omega'_{\vec{\Gamma}'}=\varepsilon_\alpha(I,J){\vert k_I\vert !\vert k_J\vert !\over
\left\vert k_{\{1,...,n\}}\right\vert !}
\omega'_{\vec{\Gamma}_1}\wedge\omega'_{\vec{\Gamma}_2}$$
et la face $\vec{F}$ (en tenant compte de son orientation)~:

$$\eqalign{\int_{\vec{F}}\omega'_{\vec{\Gamma}'}&=\varepsilon_\alpha(I,J)(-
1)^{lm_1+l+m_1}{\vert k_I\vert !\vert k_J\vert !\over \left\vert k_{\{1,...,n\}}\right\vert !}
\int_{C^+_{\{p_{i_1},...,p_{i_{n_1}}\};\{q_{l+1},...,q_{l+m_1}\}}} \omega'_{\vec{\Gamma}_1}
\cr &\hskip 4cm\int_{C^+_{A\setminus S,\ B\setminus S'\sqcup\{q\}}}
\omega'_{\vec{\Gamma}_2}.}$$
Rappelons que $\vert k_I\vert=2n_1+m_1-2$, $\vert k_J\vert=2(n-n_1)+(m-m_1+1)-2$, le signe
devant l'int\'egrale est donc~:

$$\eqalign{\varepsilon_\alpha(I,J)(-1)^{lm_1+l+m_1}&=\varepsilon_\alpha(J,I)(-1)^{|k_I||k_J|}(-1)^{lm_1+m_1+l}\cr
&=\varepsilon_\alpha(J,I)(-1)^{\vert k_I\vert\vert k_J\vert}(-1)^{lm_1+l}(-1)^{\vert k_I\vert}.}$$
Donc~:

$$c_{\vec{\Gamma}',F}=\int_{\vec{F}}\omega'_{\vec{\Gamma}'}.$$

\vskip 0.5cm
\noindent
{\bf Cas 2}~: Avec nos notations, le coefficient $c_{\vec{\Gamma}',F}$ est nul si $\vec{\Gamma}'$ ne
contient pas la fl\`eche $\overrightarrow{p_ip_j}$ et sinon~:

$$\eqalign{c_{\vec{\Gamma}',F}&=-\varepsilon_\alpha(i,j,1,...,\widehat{i,j},...,n) {(\vert
k_{\{1,...,n\}}\vert -1)! \over \left\vert k_{\{1,...,n\}}\right\vert !} \cr
&\hskip 2cm \int_{C_S}
\omega'_{\vec{\Gamma}_1} \int_{C^+{A\setminus S\sqcup\{p\},\ B}} \omega'_{\vec{\Gamma}_2}.}$$

D'autre part, on a pour la forme

$$\omega'_{\vec{\Gamma}'}=\varepsilon_\alpha(i,j,1,...,\widehat{i,j},...,n){(\vert
k_{\{1,...,n\}}\vert-1) !\over \left\vert k_{\{1,...,n\}}\right\vert !}
\omega'_{\vec{\Gamma}_1}\wedge\omega'_{\vec{\Gamma}_2}$$
et la face $\vec{F}$ (en tenant compte de son orientation)~:

$$\int_{\vec{F}}\omega'_{\vec{\Gamma}'}=
-\varepsilon_\alpha(i,j,1,...,\widehat{i,j},...,n){(\vert k_{\{1,...,n\}}\vert -1)! \over \left\vert
k_{\{1,...,n\}}\right\vert !} \int_{C_{S}} \omega'_{\vec{\Gamma}_1} \int_{C^+_{A\setminus S\sqcup \{p\},\ B}} \omega'_{\vec{\Gamma}_2},$$
avec $S=\{p_i,p_j\}$. On a donc pour tout $\vec{\Gamma}'$ et $F$~:

$$c_{\vec{\Gamma}',F}=\int_{\vec{F}}\omega'_{\vec{\Gamma}'}.$$

\vskip 0.5cm
{\bf Cas 3}~: Il n'y a pas de termes dans notre \'equation de formalit\'e dans ce cas, ou
$c_{\vec{\Gamma}',F}=0$. Mais dans ce cas, on a le lemme suivant de Kontsevich \cite {K1 \S\ 6.6.1}, \cite {Kh}~: 

$$\int_{C_S}\omega'_{\vec{\Gamma}_1}=0$$
si $|S|\ge 3$. On a donc de nouveau~:

$$c_{\vec{\Gamma}',F}=\int_{\vec{F}}\omega'_{\vec{\Gamma}'}.$$

Et ceci finit la preuve de la validit\'e de l'\'equation de formalit\'e.
\qed
\vskip 1cm
\centerline{\bf *\hskip 1cm *}
\vskip 0.5cm
\centerline{\bf *}
\paragraphe{Appendice. Formalit\'e et quantification par d\'eformation}
\qquad Nous expliquons dans ce paragraphe pourquoi la formalit\'e de Kontsevich permet d'obtenir un \'etoile-produit \`a partir d'un $2$-tenseur de Poisson. On consid\`ere une (limite projective d') alg\`ebre(s) nilpotente(s) de dimension finie $\g m$. Par exemple~:
$$\g m=\hbar\R[[\hbar]]=\lim_{\leftarrow} \hbar\R[[\hbar]]/\hbar^k\R[[\hbar]].$$
\alinea{A.1.1 Construction d'\'etoile-produits}
\qquad
On se donne une $L_\infty$-alg\`ebre $(\g g,Q)$ sur un corps $k$ de caract\'eristique z\'ero, que l'on voit comme une $Q$-vari\'et\'e formelle gradu\'ee point\'ee. Un {\sl $\g m$-point\/} de la vari\'et\'e formelle $\g g$ est par d\'efinition un morphisme de cog\`ebres~:
$$p:\g m^*\longrightarrow \Cal C(\g g).$$
Le produit tensoriel (compl\'et\'e dans le cas d'une limite projective) $\Cal C(\g g)\widehat\otimes \g m$, muni de la comultiplication de $\Cal C(\g g)$ \'etendue par $\g m$-lin\'earit\'e, admet une structure de cog\`ebre (sans co-unit\'e) sur $\g m$. On peut alors voir un $\g m$-point comme un \'el\'ement non nul de type groupe de cette cog\`ebre, c'est-\`a-dire un \'el\'ement $p\in\Cal C(\g g)\widehat\otimes \g m$
v\'erifiant~: $\Delta p=p\otimes p$. 
\prop{A.1}
Les $\g m$-points de la vari\'et\'e formelle $\g g$ sont donn\'es par~:
$$p_v=e^v-1=v+{v^2\over 2}+\cdots$$
o\`u $v$ est un \'el\'ement pair de $\g g[1]\widehat\otimes \g m$
\dem
La s\'erie a bien un sens dans  $\Cal C(\g g)\widehat\otimes \g m$. Si $p$ est un $\g m$-point, $p$ est forc\'ement pair, et on voit que la s\'erie~:
$$v=\mop{Log}(1+p)=p-{p^2\over 2}+\cdots$$
a un sens dans $\Cal C(\g g)\widehat\otimes \g m$ et d\'efinit un \'el\'ement primitif (et pair), c'est-\`a-dire que l'on a~: $\Delta v=0$. Pour d\'emontrer ce point on rajoute formellement la co-unit\'e en consid\'erant la cog\`ebre~:
$$\overline\Cal C_{\sg m}=(k.1\oplus \g m)\oplus \Cal C(\g g)\widehat\otimes \g m.$$
Un \'el\'ement de type groupe de cette cog\`ebre s'\'ecrit toujours~:
$$g=1+p$$
o\`u $p$ est de type groupe dans la cog\`ebre sans co-unit\'e $\Cal C(\g g)\widehat\otimes \g m$. Il s'agit alors de montrer que le logarithme $v$ d'un tel \'el\'ement est primitif, c'est-\`a-dire que l'on a dans $\overline\Cal C_{\sg m}$~:
$$\Delta v=v\otimes 1+1\otimes v$$
Pour cela on remarque que la cog\`ebre $\overline\Cal C_{\sg m}$ est en fait une big\`ebre. Le calcul formel suivant a alors un sens~:
$$\eqalign{\Delta\mop{Log}g     &=\mop{Log}(\Delta g)   \cr
                                &=\mop{Log}(g\otimes g) \cr
                                &=\mop{Log}\bigl((g\otimes 1)(1\otimes g)\bigr)\cr
                                &=\mop{Log}(g\otimes 1)+\mop{Log}(1\otimes g)\cr
                                &=\mop{Log}g\otimes 1+1\otimes\mop{Log}g\cr}$$

Donc, forc\'ement $v$ appartient \`a $\g g\widehat\otimes \g m$, et il est clair que $p=p_v$.
\qed
Le champ de vecteurs $Q$ s'\'etend de mani\`ere naturelle \`a $\Cal C(\g g)\widehat\otimes \g m$. Supposons que $Q$ s'annule au point $p_v$~:
$$Q(e^v-1)=0$$
On traduit ceci par le fait que $v$ v\'erifie l'{\sl \'equation de Maurer-Cartan g\'en\'eralis\'ee \/}~:
$$Q_1(v)+{1\over 2}Q_2(v.v)+\cdots =0.\leqno{\hbox{(MCG)}}$$
Si $\g g$ est une alg\`ebre de Lie diff\'erentielle gradu\'ee, \c ca se r\'eduit \`a l'\'equation de Maurer-Cartan~:
$$dv-{1\over 2}[v,v]=0.$$
(En effet $v$ est pair dans $\g g[1]\otimes \g m$, donc impair dans $\g g\otimes \g m$). Si maintenant $\Cal F$ est un $L_\infty$-morphisme entre $(\g g_1,Q)$ et $(\g g_2,Q')$, et si $v\in\g g\otimes \g m$ est tel que $Q(p_v)=0$, il est clair que~:
$$Q'\bigl(\F(p_v)\bigr)=\F\bigl(Q(p_v)\bigr)=0.$$
Or $\F(p_v)=e^w-1$ avec $w\in\g g_2\widehat\otimes \g m$ d'apr\`es la proposition A.1, puisque $\F(p_v)$ est de type groupe. Il est clair que $w$ est la projection canonique de $\F(p_v)$ sur $\g g_2\widehat\otimes m$, soit~:
$$\eqalign{w    &=\sum_{n\ge 1}{1\over n!}\F_n(v^n).\cr
                }$$
En r\'esum\'e, si $v\in\g g_1\widehat\otimes \g m$ v\'erifie (MCG), alors l'\'el\'ement $w\in \g g_2\widehat\otimes \g m$ donn\'e par l'\'egalit\'e ci-dessus v\'erifie (MCG). Dans le cas o\`u les deux $L_\infty$-alg\`ebres sont les alg\`ebres de Lie diff\'erentielles gradu\'ees des multichamps de vecteurs et des op\'erateurs polydiff\'erentiels, tout $2$-tenseur de Poisson formel~:
$$v=\hbar \gamma_1+\hbar^2 \gamma_2+\cdots$$
donne naissance gr\^ace \`a ce processus \`a un op\'erateur bidiff\'erentiel formel $w$ tel que $\mu+w$ soit un \'etoile-produit, $\mu$ d\'esignant la multiplication usuelle de deux fonctions.
\alinea{A.2. Equivalence des foncteurs de d\'eformation}
\qquad On suppose toujours le corps de base $k$ de caract\'eristique z\'ero. Soit $\g g$ une alg\`ebre de Lie diff\'erentielle gradu\'ee. Rappelons \cite {K1 \S~3.2} que le foncteur de d\'eformation $\mop{Def}_{\sg g}$ associe \`a toute alg\`ebre commutative nilpotente de dimension finie $\g m$ l'ensemble des classes de solutions de degr\'e $1$ de l'\'equation de Maurer-Cartan dans $\g g\otimes \g m$ modulo l'action du groupe de jauge, c'est-\`a-dire le groupe nilpotent $G_{\sg m}=\exp (\g g^0\otimes \g m)$, dont l'action (par des transformations affines de l'espace $\g g^1 \otimes \g m$) est donn\'ee infinit\'esimalement par~:
$$\alpha.\gamma=d\alpha+[\alpha,\gamma]$$
pour tout $\alpha\in\g g^0\otimes\g m$ et pour tout $\gamma\in\g g^1\otimes \g m$. Ce foncteur s'\'etend naturellement aux limites projectives d'alg\`ebres commutatives nilpotentes de dimension finie~: $\mop{Def}_{\sg g}(\g m)$ est dans ce cas d\'efini comme l'ensemble des classes de solutions de degr\'e $1$ de l'\'equation de Maurer-Cartan dans le produit tensoriel compl\'et\'e $\g g\widehat\otimes \g m$ modulo l'action du groupe pro-nilpotent $G_{\sg m}=\exp (\g g^0\widehat\otimes \g m)$.
\ssq
L'\'equivalence de jauge peut aussi se d\'efinir pour une $L_\infty$-alg\`ebre quelconque $(\g g,Q)$~: deux solutions de l'\'equation de Maurer-Cartan g\'en\'eralis\'ee $\gamma_0$ et $\gamma_1$ dans $\g g^1\otimes \g m$ sont \'equivalentes s'il existe une famille polynomiale $\xi(t)_{t\in k}$ de champs de vecteurs de degr\'e $-1$ et une famille polynomiale $\gamma(t)_{t\in k}$ de solutions de l'\'equation de Maurer-Cartan g\'en\'eralis\'ee dans $\g g^1\otimes \g m$ telles que~:
$$\eqalign{{d\gamma(t)\over dt}      &=[Q,\xi(t)]\bigl(\gamma (t)\bigr)\cr
             \gamma(0)     &=\gamma_0,\hbox to 12mm{}\gamma(1)=\gamma_1.\cr}\eqno{(**)}$$
On v\'erifie facilement que cette relation est une relation d'\'equivalence, ce qui permet de d\'efinir le foncteur de d\'eformation $\mop{Def}_{\sg g}$ comme la correspondance qui \`a toute alg\`ebre commutative nilpotente de dimension finie $\g m$ associe l'ensmble $\mop{Def}_{\sg g}(\g m)$ des classes de solutions de degr\'e $1$ de l'\'equation de Maurer-Cartan g\'en\'eralis\'ee dans $\g g\otimes \g m$ modulo l'\'equivalence de jauge.
\prop{A.2.1 \rm(cf. \cite {K1 \S~4.5.2})}
1). Dans le cas d'une alg\`ebre de Lie diff\'erentielle gradu\'ee les deux notions d'\'equivalence de jauge (et donc de foncteur de d\'eformation) co\"\i ncident.
\smallskip
2). Soient $\g g_1$ et $\g g_2$ deux alg\`ebres de Lie diff\'erentielles gradu\'ees. Alors le foncteur $\mop{Def}_{\sg g_1\oplus\sg g_2}$ est naturellement \'equivalent au produit des foncteurs $\mop{Def}_{\sg g_1}\times \mop{Def}_{\sg g_2}$.
\smallskip
3). Le foncteur de d\'eformation est trivial pour une $L_\infty$-alg\`ebre lin\'eaire contractile.
\dem
Les points 2) et 3) sont faciles \`a \'etablir. Pour \'etablir le premier point on remarque que l'action de $\alpha\in\g g^0\otimes \g m$ sur $\g g^1\otimes \g m$ est donn\'ee par le champ de vecteurs~:
$$D_\alpha=[Q,R_\alpha],$$
o\`u $R_\alpha$ est le champ de vecteurs constant \'egal \`a $\alpha$. Ce champ de vecteurs est bien de degr\'e $-1$. L'\'equivalence de jauge au sens des alg\`ebres de Lie diff\'erentielles gradu\'ees $\gamma_1=(\exp \alpha).\gamma_0$ entra\^\i ne donc l'\'equivalence de jauge au sens des $L_\infty$-alg\`ebres, avec $\xi(t)=R_\alpha$ pour tout $t$ et $\gamma(t)=(\exp t\alpha).\gamma_0$.
\ssq
Supposons maintenant que $\gamma_0$ et $\gamma_1$ sont \'equivalents au sens des $L_\infty$-alg\`ebres. Soient $\xi(t)$ et $\gamma(t)$ les familles polynomiales de champs de vecteurs de degr\'e $-1$ et de solutions de l'\'equation de Maurer-Cartan g\'en\'eralis\'ee respectivement, telles que l'\'equation (**) soit v\'erifi\'ee. Le champ de vecteurs de degr\'e z\'ero $[Q,\xi(t)]$ s'\'ecrit explicitement en tout $\g m$-point $\gamma\in\g g^1\otimes \g m$~:
$$[Q,\xi(t)](\gamma)=d\xi(t)+[\xi(t),\gamma].$$
On effectue ce calcul en appliquant la cod\'erivation $[Q,\xi(t)]$ \`a l'\'el\'ement de type groupe $e^\gamma -1$. Compte tenu de l'\'equation (**) on voit que le vecteur tangent ${d \gamma(t)\over dt}$ au point $\gamma(t)$ est donn\'e par un champ de vecteurs provenant de l'action d'un \'el\'ement de l'alg\`ebre de Lie $\g g^0\otimes \g m$.
\ssq
Montrons par r\'ecurrence sur le degr\'e $d$ de $\gamma_t$ (en tant que polyn\^ome) qu'il existe un entier $r$ (d\'ependant de $d$) tel que pour tout $t\in k$ il existe $g_t=\exp (tD_1+\cdots +t^r D_r)\in G_{\sg m}$ tel que $\gamma_t=g_t.\gamma_0$~: si $d=0$ on a $\gamma_t=\gamma_0$ et la constante $g_t=\mop{Id}$ convient. Supposons donc que la propri\'et\'e soit vraie au rang $d-1$. Supposons que $\gamma_t$ soit un polyn\^ome de degr\'e $d$, que l'on peut \'ecrire~:
$$\gamma_t=\tilde \gamma_t+t^d\gamma_d.$$
Gr\^ace \`a l'hypoth\`ese de r\'ecurrence on peut \'ecrire~:
$$\eqalign{\gamma_t    &=e^{t^d\gamma_d}\tilde \gamma_t     \cr
&=e^{t^d\gamma_d}e^{tD_1+\cdots +t^r D_r}\gamma_0      .\cr}$$
Le terme $\gamma_d$ appartient \`a l'alg\`ebre de Lie $\g g^0\otimes \g m$. La s\'erie de Campbell-Hausdorff ne comprend qu'un nombre fini de termes dans le cas d'un groupe nilpotent, ce qui permet de conclure. Le passage aux limites projectives d'alg\`ebres commutatives nilpotentes de dimension finie se fait sans difficult\'e.
\qed
Compte tenu de la proposition $V.2$, la proposition A.2.1 entra\^\i ne le r\'esultat suivant~:
\th{A.2.2}
Soient $\g g_1$ et $\g g_2$ deux $L_\infty$-alg\`ebres quasi-isomorphes. Alors les foncteurs de d\'eformation de $\g g_1$ et de $\g g_2$ sont isomorphes.
\ndem
En particulier les classes d'\'equivalence de jauge de $2$-tenseurs de Poisson formels sur une vari\'et\'e sont en bijection avec les classes d'\'equivalence de jauge d'\'etoile-produits.
\paragraphe{R\'ef\'erences~:}
\bib{AKSZ}M. Alexandrov, M. Kontsevich, A. Schwarz, O. Zaboronsky, {\sl The geometry of the Master equation and topological quantum field theory\/}, Int. J. Mod. Phys. et hep-th 9502010.
\bib{BFFLS}F. Bayen, M. Flato, C. Fr\o nsdal, A. Lichnerowicz, D. Sternheimer, {\sl Deformation theory and quantization I. Deformations of symplectic structures}, Ann. Phys. 111 No1, 61-110 (1978).
\bib{FM}W. Fulton, R. MacPherson, {\sl Compactification of configuration spaces}, Ann. Math. 139, 183-225 (1994).  
\bib {H-S}V. Hinich, V. Schechtman, {\sl Homotopy Lie algebras\/}, I.M. Gelfand Seminar, Adv. Sov. math. 16 (2), 1993.
\bib {Kh}A.G. Khovanski\u\i, {\sl On a lemma of Kontsevich.\/} Funct.
Anal. Appl. 31, no. 4, 296--298 (1998) 
\bib{K1}M. Kontsevich, {\sl Deformation quantization of Poisson manifolds\/}, q-alg. 9709040.
\bib{K2}M. Kontsevich, {\sl Formality conjecture}, D. Sternheimer et al. (eds.), Deformation theory and symplectic geometry, Kluwer, 139-156 (1997).
\bib MS. Majid, {\sl Algebras and Hopf algebras in braided categories\/}, q-alg. 9509023.
\bib QD. Quillen, {Rational homotopy theory\/}, Ann. Math. 90 (1969), 205-295.
 
\bye

\bye